\documentclass[11pt]{article}
\usepackage{latexsym}
\usepackage{amsmath,amssymb}

\usepackage{graphicx}

\begin{document}

\title{The analogue of Wilczynski's projective frame in Lie sphere geometry:
Lie-applicable surfaces and commuting Schr\"odinger operators
with magnetic fields.}
\author{{\Large Ferapontov E.V. \ }  \\
    Department of Mathematical Sciences \\
    Loughborough University \\
    Loughborough, Leicestershire LE11 3TU \\
    United Kingdom \\
    e-mail: {\tt E.V.Ferapontov@lboro.ac.uk} \\
    and \\
    Centre for Nonlinear Studies \\
    Landau Institute of Theoretical Physics\\
    Academy of Science of Russia, Kosygina 2\\
    117940 Moscow, GSP-1, Russia\\
    e-mail: {\tt fer@landau.ac.ru}\\
}
\date{}
\maketitle

\newtheorem{theorem}{Theorem}

\pagestyle{plain}

\maketitle

\begin{abstract}

We propose a twistor construction of surfaces in Lie sphere geometry
based on the linear system which copies equations of
Wilczynski's projective frame.
In the particular case of Lie-applicable surfaces this linear
system describes joint eigenfunctions of a pair of commuting Schr\"odinger
operators 
with magnetic terms.

\bigskip


\bigskip

\end{abstract}

\newpage
\section{Introduction}

Let $M^2$ be a surface in $E^3$ parametrized by  coordinates $R^1, R^2$
of curvature lines, with the radius-vector ${\bf r}$ and the unit normal
${\bf n}$ satisfying the Weingarten equations
\begin{equation}
\begin{array}{c}
\partial_1${\bf r}$ =w^1 \partial_1${\bf n}$ \\
\ \\
\partial_2${\bf r}$ =w^2 \partial_2${\bf n}$, \\
\end{array}
\label{Weingarten}
\end{equation}
where $w^1$ and $ w^2$
are the radii of principal curvature, $\partial_1=\partial/\partial R^1, \
\partial_2=\partial/\partial R^2$.
Let us recall the construction of the Lie sphere map \cite{Lie}. With any
sphere
$S(R, {\bf c})$ having radius $R$ and  center ${\bf c}=(c^1, c^2, c^3)$ this
map
associates the 6-vector
$\{ y_0, y_1, y_2, y_3, y_4, y_5\} $  with hexaspherical coordinates
$$
\left \{ \frac{1+{\bf c}^2-R^2}{2}, ~~ \frac{1-{\bf c}^2+R^2}{2}, ~ {\bf c},
~ R\right \},
$$
which obey the relation
\begin{equation}
-y_0^2+y_1^2+y_2^2+y_3^2+y_4^2-y_5^2=0.
\label{Liequadric}
\end{equation}
This equation defines the so-called Lie quadric. Thus, with any sphere
$S(R, {\bf c})$ in $E^3$ we associate a point
on the Lie quadric (\ref{Liequadric}).
The Lie sphere map linearises the action of the Lie sphere group which is a
group
of contact transformations in $E^3$ generated by conformal transformations
and normal 
shifts. In hexaspherical coordinates, the action of the Lie sphere group
coincides with the linear action of $SO(4, 2)$ which preserves the Lie
quadric
(\ref{Liequadric}).
The reader may consult \cite{Lie},
\cite{Blaschke}, \cite{Cecil}, \cite{Pinkall1}, \cite{Pinkall2}
for further properties of this construction.
Applying  Lie sphere map to the curvature spheres
$S(w^1, {\bf r}-w^1{\bf n})$ and  $S(w^2, {\bf r}-w^2{\bf n})$ of the
surface
$M^2$, we obtain a pair of two-dimensional submanifolds
$$
{\bf U}=\left \{\frac{1+{\bf r}^2-2w^1({\bf r}, {\bf n})}{2}, ~~
\frac{1-{\bf r}^2+2w^1({\bf r}, {\bf n})}{2}, ~~ {\bf r}-w^1{\bf n}, ~~ w^1
\right \}
$$
and
$$
{\bf V}=\left \{\frac{1+{\bf r}^2-2w^2({\bf r}, {\bf n})}{2}, ~~
\frac{1-{\bf r}^2+2w^2({\bf r}, {\bf n})}{2}, ~~ {\bf r}-w^2{\bf n}, ~~ w^2
\right \}
$$
of the Lie quadric. Blaschke's approach  \cite{Blaschke} to the Lie sphere
geometry of surfaces
in 3-space is based on the following two simple facts:

\noindent (a) By  construction, vectors ${\bf U}$ and ${\bf V}$ have zero
norm
$$
({\bf U}, \  {\bf U})= ({\bf V}, \  {\bf V})=0
$$
in the scalar product defined by (\ref{Liequadric}).

\noindent (b) The triple
${\bf U}, \ \partial_2 {\bf U}, \  \partial_2^2{\bf U}$
is orthogonal to the triple
${\bf V}, \ \partial_1 {\bf V}, \  \partial_1^2{\bf V}$.

\noindent Using appropriate linear combinations among these triples, one can
construct an invariant 6-frame canonically associated with a surface $M^2$
(Appendix B, see also \cite{Blaschke}, \cite{Fer001}).
Conversely, given a pair of 6-vectors ${\bf U}$
and ${\bf V}$ satisfying (a) and (b),
the surface can be reconstructed uniquely
as the envelope of the corresponding family of curvature spheres.

Our construction of vectors ${\bf U}$ and ${\bf V}$ which
satisfy the above properties
was borrowed  from projective differential geometry.
In 1907, Wilczynski \cite{Wilczynski} proposed the approach to surfaces in
projective space 
based on a linear system
\begin{equation}
\begin{array}{c}
{\bf r}_{xx}=\beta \ {\bf r}_y+\frac{1}{2}(V-\beta_y) \ {\bf r} \\
\ \\
{\bf r}_{yy}=\gamma \ {\bf r}_x+\frac{1}{2}(W-\gamma_x) \ {\bf r}
\end{array}
\label{r}
\end{equation}
where $\beta, \gamma, V, W$ are real functions of $x$ and  $y$.
Cross-differentiating (\ref{r}) and assuming
${\bf r}, {\bf r}_x, {\bf r}_y, {\bf r}_{xy}$ to be independent, we arrive
at
the compatibility conditions
\begin{equation}
\begin{array}{c}
\beta_{yyy}-2\beta_yW-\beta W_y=
\gamma_{xxx}-2\gamma_xV-\gamma V_x\\
\ \\
W_x=2\gamma \beta_y+\beta \gamma_y \\
\ \\
V_y=2\beta \gamma_x+\gamma \beta_x.
\end{array}
\label{GC1}
\end{equation}
For any
$\beta, \gamma, V, W$ satisfying (\ref{GC1}) the
linear system (\ref{r}) has rank 4, so that
${\bf r}=(r^0:r^1:r^2:r^3)$ can be viewed as a 4-component real vector of
homogeneous coordinates of a
surface $M^2$ in a projective space ${\bf RP}^3$. The compatibility
conditions (\ref{GC1})
can be interpreted as the `Gauss-Codazzi' equations in projective
differential
 geometry.
The independent variables $x, y$ play the role of asymptotic coordinates on
the surface $M^2$. Introducing the vectors
$$
{\cal U}={\bf r}\wedge{\bf r}_x ~~ {\rm and} ~~
{\cal V}={\bf r}\wedge{\bf r}_y
$$
in $\Lambda^2(R^4)$, one   readily verifies that

\noindent (a) vectors ${\cal U}$ and ${\cal V}$ have zero norm
$$
({\cal U}, \  {\cal U})= ({\cal V}, \  {\cal V})=0
$$
in the natural scalar product in $\Lambda^2(R^4)$ defined by the
Pl\"ucker formula;

\noindent (b) the triple
${\cal U}, \ {\cal U}_y, \  {\cal U}_{yy}$
is orthogonal to the triple
${\cal V}, \  {\cal V}_x, \  {\cal V}_{xx}$.

\noindent The passage from a projective surface to a pair of submanifolds
${\cal U}$ and ${\cal V}$ in the Pl\"ucker quadric in $\Lambda^2(R^4)$
plays an important role in projective differential geometry.
We refer to \cite{Wilczynski}, \cite{Bol}, \cite{Lane},
\cite{Fer8}, \cite{Fer002}, \cite{Sasaki1}) for a further discussion.
Some more details on  Wilczynski's
approach are included in  Appendix A.

\bigskip

The main observation of this paper is that a similar approach,
 based on the linear system
\begin{equation}
\begin{array}{c}
\partial_1^2{ {\mbox{\boldmath $\psi$}}}=
-i p \ \partial_2{{\mbox{\boldmath $\psi$}}}+\frac{1}{2}
(V+i\partial_2p) \ {\mbox{\boldmath $\psi$}} \\
\ \\
\partial_2^2{\mbox{\boldmath $\psi$}}=
i q \ \partial_1{\mbox{\boldmath $\psi$}}+\frac{1}{2}
(W-i\partial_1q) \ {\mbox{\boldmath $\psi$}},
\end{array}
\label{LieWil}
\end{equation}
applies to Lie sphere geometry. Here $\partial_1=\partial/\partial R^1$ and
$ \partial_2=\partial/\partial R^2$
denote 
partial derivatives with respect to the independent variables $R^1, \ R^2$,
and $p, q, V, W$ are real potentials.
This system  is a complex analogue of the linear system (\ref{r})
(indeed, the complex transformation
$\beta=-ip, \ \gamma = iq$ identifies both systems).
Again,  cross-differentiation produces the compatibility conditions
\begin{equation}
\begin{array}{c}
\partial_2^3p-2W\partial_2p-p\partial_2 W+
\partial_1^3q-2V\partial_1q-q\partial_1 V=0\\
\ \\
\partial_1W=2q \partial_2p+p \partial_2q \\
\ \\
\partial_2V=2p \partial_1q+q \partial_1p.
\end{array}
\label{VW}
\end{equation}
For any fixed
$p, q, V, W$ satisfying (\ref{VW}), the linear system
(\ref{LieWil})  is compatible
and its  solution space has complex dimension 4, so that we can view
${\mbox{\boldmath $\psi$}}$
as an element of the twistor space ${\bf CP}^3$.
In what follows, equations (\ref{VW}) will be identified with
the `Gauss-Codazzi' equations in Lie sphere geometry, while the independent
variables
$R^1, R^2$ will play the role of  curvature line coordinates.
An important property of linear system (\ref{LieWil}) is the existence of
the invariant 
pseudo-Hermitian scalar product $(\ , \ )$ of the signature $(2, 2)$ such
that
\begin{equation}
\begin{array}{c}
({\mbox{\boldmath $\psi$}},\ \partial_1\partial_2 {\mbox{\boldmath
$\psi$}})=
(\partial_1\partial_2 {{\mbox{\boldmath $\psi$}}}, \ {{\mbox{\boldmath
$\psi$}}})=-1,
~~~~
(\partial_1 {\mbox{\boldmath $\psi$}},\ \partial_2 {\mbox{\boldmath
$\psi$}})=
(\partial_2 {\mbox{\boldmath $\psi$}}, \ \partial_1 {\mbox{\boldmath
$\psi$}})=1 \\
 \ \\
(\partial_1\partial_2 {\mbox{\boldmath $\psi$}},\
\partial_1\partial_2 {\mbox{\boldmath $\psi$}})=-pq
\end{array}
\label{norm}
\end{equation}
(all other scalar products being zero). Fixing ${\mbox{\boldmath $\psi$}}\in
{\bf C}^4$ 
satisfying both
(\ref{LieWil}) and (\ref{norm}), we define two real vectors
${\bf U}$ and ${\bf V}$ in $\Lambda^2({\bf C}^4)$ as
$$
{\bf U}=-2\ Im ({\mbox{\boldmath $\psi$}}\wedge
\partial_1{{\mbox{\boldmath $\psi$}}}), ~~~
{\bf V}=2\ Re ({{\mbox{\boldmath $\psi$}}}\wedge
\partial_2{{\mbox{\boldmath $\psi$}}}).
$$
Introducing  in $\Lambda^2({\bf C}^4)$ a pseudo-Hermitian scalar product $(\
, \ )$
induced by (\ref{norm}) (the wedge product $\wedge$ and the scalar product
$(\ , \ )$  in $\Lambda^2({\bf C}^4)$ are explicitly defined in Appendix C),
we show that 
the restriction of $(\ , \ )$ to the 6-dimensional invariant real subspace
in
$\Lambda^2({\bf C}^4)$ spanned by ${\bf U}, \ \partial_2 {\bf U}, \
\partial_2^2 {\bf U}$ and ${\bf V}, \ \partial_1 {\bf V}, \
\partial_1^2 {\bf V}$ is a real scalar product of the signature $(4, 2)$.
Moreover,

\noindent (a) vectors ${\bf U}$ and ${\bf V}$ have zero norm
$$
({\bf U}, \  {\bf U})= ({\bf V}, \  {\bf V})=0;
$$

\noindent (b) the triple
${\bf U}, \ \partial_2 {\bf U}, \  \partial_2^2{\bf U}$
is orthogonal to the triple
${\bf V}, \ \partial_1 {\bf V}, \  \partial_1^2{\bf V}$.

\noindent Thus, ${\bf U}$ and  ${\bf V}$ are hexaspherical coordinates of
curvature spheres of a surface $M^2$ parametrized by curvature line
coordinates
$R^1, R^2$. In fact, our construction is based on the
isomorphism  
$SU(2, 2)/{\pm 1} \to SO(4, 2)$ which is a basis of twistor theory
\cite{Penrose}, \cite{Tod}, \cite{Mason}.

The main motivation for the construction described above
comes from  the study of Lie-applicable surfaces. We recall that two
surfaces 
are called Lie-applicable if,
being non-equivalent under Lie sphere transformations,
they have the same coefficients $p$ and $q$ in the appropriate curvature
line parametrization $R^1, R^2$ \cite{Blaschke}. Analytically,
Lie-applicable surfaces
are described by equations (\ref{VW}) which, for given $p$ and $q$, are not
uniquely solvable
for $V$ and $W$. Examples presented in section 3 demonstrate that for some
particularly interesting classes of Lie-applicable surfaces, equations
(\ref{LieWil}) describe joint eigenfunctions of a pair of commuting
Schr\"odinger operators with magnetic fields. These examples include
nonsingular
doubly periodic operators on a two-dimensional torus and the Dirac monopole
on a 
two-dimensional sphere \cite{Wu}.
It should be emphasized that
this important structure is  not visible in
the standard approach based on  Blaschke's 6-frame,  becoming apparent only
after applying the twistor construction.

In section 4 we give a separate treatment of canal surfaces
(which are characterized by a condition $q=0$ (or $p=0$)),
since the construction of the Lie sphere frame
adopted in section 2 does not automatically carry over to this case. As an
example, we 
explicitly construct the surfaces of revolution for which equations
(\ref{LieWil})
reduce to  eigenfunction equations of the Schr\"odinger operator
in a homogeneous magnetic field. The case of the first nontrivial
 Landau level is worked out in detail.

\section{Analog of Wilczynski's projective frame in Lie sphere geometry}

In this section we propose an approach to surfaces in Lie sphere
geometry based on the linear system (\ref{LieWil}) satisfying the
compatibility conditions
(\ref{VW}).
Our further constructions will follow those from projective
differential geometry --- see Appendix A.
Notice first that system (\ref{LieWil}) is covariant under
transformations of the form
\begin{equation}
(R^1)^*=f(R^1), ~~~~ (R^2)^*=g(R^2), ~~~~
{{\mbox{\boldmath $\psi$}}}^*=\sqrt {f'(R^1)g'(R^2)}~{{\mbox{\boldmath
$\psi$}}}
\label{newR}
\end{equation}
which act on the potentials $p, q, V, W$ as follows
\begin{equation}
\begin{array}{c}
p^{*}=p g'/(f')^2, ~~~~  V^{*}(f')^2=V+S(f)\\
\ \\
q^{*}=q f'/(g')^2, ~~~~ W^{*}(g')^2=W+S(g).
\end{array}
\label{newpq}
\end{equation}
Here $S(\,\cdot\,)$ is the  Schwarzian derivative
$$
S(f)=\frac{f'''}{f'} - \frac{3}{2} \left(\frac{f''}{f'}\right)^2,
$$
so that $V$ and $W$ can be interpreted as projective connections.
Transformation formulae (\ref{newR}), (\ref{newpq}) imply that
the symmetric 2-form
\begin{equation}
-pq \ dR^1dR^2
\label{metric}
\end{equation}
and the conformal class of the cubic form
\begin{equation}
p(dR^1)^3-q(dR^2)^3
\label{cubic}
\end{equation}
are invariant. In what follows they will play the roles
of the Lie-invariant metric and the Lie-invariant cubic form, respectively.
In this section we assume both $p$ and  $q$ to be  nonzero.
Let us introduce the four vectors
\begin{equation}
\begin{array}{c}
{{\mbox{\boldmath $\psi$}}}, ~~~
{{\mbox{\boldmath $\psi$}}}_1=\partial_1{{\mbox{\boldmath
$\psi$}}}-\frac{1}{2}\frac{\partial_1q}{q}{{\mbox{\boldmath $\psi$}}}, ~~~
{{\mbox{\boldmath $\psi$}}}_2=\partial_2{{\mbox{\boldmath
$\psi$}}}-\frac{1}{2}\frac{\partial_2p}{p}{{\mbox{\boldmath $\psi$}}}, \\
\ \\
{\mbox{\boldmath $\eta$}}=\partial_1\partial_2{{\mbox{\boldmath $\psi$}}}-
\frac{1}{2}\frac{\partial_1q}{q}\partial_2{{\mbox{\boldmath $\psi$}}}-
\frac{1}{2}\frac{\partial_2p}{p}\partial_1{{\mbox{\boldmath $\psi$}}}+
\left(\frac{1}{4}\frac{\partial_2p\partial_1q}{pq} -
\frac{1}{2}{pq}\right){{\mbox{\boldmath $\psi$}}} \\
\end{array}
\label{framelie}
\end{equation}
which are straightforward analogues of vertices of
Wilczynski's moving tetrahedral --- see formulae (\ref{frame}) in  Appendix
A.
Notice that under  transformations
(\ref{newR}) vectors (\ref{framelie}) acquire  nonzero multiples which do
not
change them as points in the complex projective space.
Using (\ref{LieWil}) and (\ref{framelie}),
we readily derive for
${{\mbox{\boldmath $\psi$}}}, {{\mbox{\boldmath $\psi$}}}_1,
 {{\mbox{\boldmath $\psi$}}}_2, {\mbox{\boldmath $\eta$}}$
the linear system
\begin{equation}
\begin{array}{c}
\partial_1\left(\begin{array}{c}
{{\mbox{\boldmath $\psi$}}}\\
{{\mbox{\boldmath $\psi$}}}_1\\
{{\mbox{\boldmath $\psi$}}}_2\\
{\mbox{\boldmath $\eta$}}
\end{array}\right)=
\left(\begin{array}{cccc}
\frac{1}{2}\frac{\partial_1q}{q}& 1&0&0\\
\frac{1}{2}b & -\frac{1}{2}\frac{\partial_1q}{q}& -ip &0\\
\frac{1}{2}k&0&\frac{1}{2}\frac{\partial_1q}{q}&1\\
-\frac{i}{2}p
a&\frac{1}{2}k&\frac{1}{2}b&-\frac{1}{2}\frac{\partial_1q}{q}
\end{array}\right)
\left(\begin{array}{c}
{{\mbox{\boldmath $\psi$}}}\\
{{\mbox{\boldmath $\psi$}}}_1\\
{{\mbox{\boldmath $\psi$}}}_2\\
{\mbox{\boldmath $\eta$}}
\end{array}\right)\\
\ \\
\partial_2\left(\begin{array}{c}
{{\mbox{\boldmath $\psi$}}}\\
{{\mbox{\boldmath $\psi$}}}_1\\
{{\mbox{\boldmath $\psi$}}}_2\\
{\mbox{\boldmath $\eta$}}
\end{array}\right)=
\left(\begin{array}{cccc}
\frac{1}{2}\frac{\partial_2p}{p}& 0&1&0\\
\frac{1}{2}l & \frac{1}{2}\frac{\partial_2p}{p}& 0&1\\
\frac{1}{2}a&iq &-\frac{1}{2}\frac{\partial_2p}{p}&0\\
\frac{i}{2}q b&\frac{1}{2}a
&\frac{1}{2}l&-\frac{1}{2}\frac{\partial_2p}{p}
\end{array}\right)\left(
\begin{array}{c}
{{\mbox{\boldmath $\psi$}}}\\
{{\mbox{\boldmath $\psi$}}}_1\\
{{\mbox{\boldmath $\psi$}}}_2\\
{\mbox{\boldmath $\eta$}}
\end{array}\right)
\end{array}
\label{lieframe}
\end{equation}
where we introduced the notation
$$
\begin{array}{c}
k=pq - \partial_1\partial_2(\ln p), ~~~~
l=pq - \partial_1\partial_2(\ln q)\\
\ \\
a=W-\partial_2^2(\ln p)-\frac{1}{2}(\partial_2 \ln p)^2, ~~~~
b=V-\partial_1^2(\ln q)-\frac{1}{2}(\partial_1 \ln q)^2.
\end{array}
$$
The compatibility conditions of equations (\ref{lieframe}) imply
\begin{equation}
\begin{array}{c}
\partial_1\partial_2 \ln p=pq -k, ~~~~
\partial_1\partial_2 \ln q=pq -l\\
\ \\
\partial_1a=\partial_2k+\frac{\partial_2p}{p}\, k, ~~~~
\partial_2b=\partial_1l+\frac{\partial_1q}{q}\, l \\
\ \\
p\,  \partial_2a+2\, a\, \partial_2p +
q\,  \partial_1b + 2\, b\, \partial_1q = 0
\end{array}
\label{LieGC}
\end{equation}
which is just an equivalent form of
the `Gauss-Codazzi' equations (\ref{VW}).
An important property of  system (\ref{lieframe}) is the existence of the
quadratic integral
\begin{equation}
-{{\mbox{\boldmath $\psi$}}}\bar {\mbox{\boldmath $\eta$}}+{{\mbox{\boldmath
$\psi$}}}_1{\bar {{\mbox{\boldmath $\psi$}}}}_2+
{{\mbox{\boldmath $\psi$}}}_2 {\bar {{\mbox{\boldmath $\psi$}}}}_1 -
{\mbox{\boldmath $\eta$}}{\bar {{\mbox{\boldmath $\psi$}}}}
\label{integral1}
\end{equation}
which defines an invariant pseudo-Hermitian scalar product of the signature
$(2, 2)$ on the space of solutions of system (\ref{lieframe}).
Using  (\ref{framelie}), this integral can be rewritten in the form
\begin{equation}
-{{\mbox{\boldmath $\psi$}}}\ \partial_1\partial_2{\bar{{\mbox{\boldmath
$\psi$}}}}+ 
\partial_1{{\mbox{\boldmath $\psi$}}}\ \partial_2 {\bar {{\mbox{\boldmath
$\psi$}}}}+
\partial_2{{\mbox{\boldmath $\psi$}}}\ \partial_1{\bar {{\mbox{\boldmath
$\psi$}}}}-
{\bar {{\mbox{\boldmath $\psi$}}}}\ \partial_1\partial_2{{{\mbox{\boldmath
$\psi$}}}}+pq\ {{\mbox{\boldmath $\psi$}}}{\bar {{\mbox{\boldmath
$\psi$}}}}.
\label{integral2}
\end{equation}
The existence of the invariant pseudo-Hermitian scalar product
(\ref{integral1})
allows to choose a basis of solutions of (\ref{lieframe})
such that 
\begin{equation}
({{\mbox{\boldmath $\psi$}}}_1, \ {{\mbox{\boldmath
$\psi$}}}_2)=({{\mbox{\boldmath $\psi$}}}_2, \ {{\mbox{\boldmath
$\psi$}}}_1)=1, ~~~~
({{\mbox{\boldmath $\psi$}}}, \ {\mbox{\boldmath $\eta$}})=({\mbox{\boldmath
$\eta$}}, 
\ {{\mbox{\boldmath $\psi$}}})=-1,
\label{Herm}
\end{equation}
all other scalar products being zero.
Notice that formulae (\ref{Herm}) are equivalent to  (\ref{norm}).
Here $(~,~)$ denotes the pseudo-Hermitian scalar product in
${\bf C}^4$
$$
({\bf a}, {\bf b})=-a^0\bar b^3+a^1\bar b^2+a^2\bar b^1-a^3\bar b^0
$$
of the signature $(2, \ 2)$. Equations (\ref{lieframe}) also imply that
the determinant
${{\mbox{\boldmath $\psi$}}}\wedge {{\mbox{\boldmath $\psi$}}}_1 \wedge
{{\mbox{\boldmath $\psi$}}}_2\wedge {\mbox{\boldmath $\eta$}}$
is invariant:
$$
\partial_1 ({{\mbox{\boldmath $\psi$}}}\wedge {{\mbox{\boldmath $\psi$}}}_1
\wedge {{\mbox{\boldmath $\psi$}}}_2\wedge {\mbox{\boldmath $\eta$}})=
\partial_2 ({{\mbox{\boldmath $\psi$}}}\wedge {{\mbox{\boldmath $\psi$}}}_1
\wedge {{\mbox{\boldmath $\psi$}}}_2\wedge {\mbox{\boldmath $\eta$}})=0
$$
(indeed, both matrices in (\ref{lieframe}) are traceless). Thus, besides
(\ref{Herm}), we can impose the additional constraint
\begin{equation}
{{\mbox{\boldmath $\psi$}}}\wedge {{\mbox{\boldmath $\psi$}}}_1 \wedge
{{\mbox{\boldmath $\psi$}}}_2\wedge {\mbox{\boldmath $\eta$}}=1.
\label{det}
\end{equation}
From now on, we fix a null-tetrad
${{\mbox{\boldmath $\psi$}}}, \ {{\mbox{\boldmath $\psi$}}}_1, \
{{\mbox{\boldmath $\psi$}}}_2, \ {\mbox{\boldmath $\eta$}}$
satisfying both (\ref{Herm}) and (\ref{det}).  Notice that this basis
${{\mbox{\boldmath $\psi$}}}, \ {{\mbox{\boldmath $\psi$}}}_1, \
{{\mbox{\boldmath $\psi$}}}_2, \ {\mbox{\boldmath $\eta$}}$
described above is defined up to a natural linear action of the group
$SU(2, 2)$ which preserves both (\ref{Herm}) and (\ref{det}).
Introducing the basis in $\Lambda^2{\bf C}^4$ as follows (notice the analogy
with
self-dual and anti-self-dual forms in twistor theory)
$$
\begin{array}{c}
{\cal U}=i \ {{\mbox{\boldmath $\psi$}}} \wedge {{\mbox{\boldmath
$\psi$}}}_1, ~~~ {\cal V}={{\mbox{\boldmath $\psi$}}} \wedge
{{\mbox{\boldmath $\psi$}}}_2, \\
\ \\
{\cal A} =i \ {{\mbox{\boldmath $\psi$}}}_2\wedge {{\mbox{\boldmath
$\psi$}}}_1+i \ {{\mbox{\boldmath $\psi$}}}\wedge {\mbox{\boldmath $\eta$}},
~~~
{\cal B} ={{\mbox{\boldmath $\psi$}}}_1\wedge {{\mbox{\boldmath
$\psi$}}}_2+{{\mbox{\boldmath $\psi$}}}\wedge {\mbox{\boldmath $\eta$}}, \\
\ \\
{\cal P} = 2i\, {{\mbox{\boldmath $\psi$}}}_2\wedge {\mbox{\boldmath
$\eta$}}, ~~~
{\cal Q} = 2\, {{\mbox{\boldmath $\psi$}}}_1\wedge {\mbox{\boldmath
$\eta$}},
\end{array}
$$
we arrive at the  equations
\begin{equation}
\begin{array}{c}
\partial_1\left(\begin{array}{c}
{\cal U}\\
{\cal A}\\
{\cal P}\\
{\cal V}\\
{\cal B}\\
{\cal Q}
\end{array}\right)=
\left(\begin{array}{cccccc}
0 & 0 & 0 & p & 0 & 0\\
k & 0 & 0 & 0 & 0 & 0\\
0 & k & 0 & -pa & 0 & 0\\
0 & 0 & 0 & \frac{\partial_1q}{q} & 1 & 0\\
0 & 0 & 0 & b & 0 & 1\\
pa & 0 & -p & 0 & b &-\frac{\partial_1q}{q}
\end{array}\right)
\left(\begin{array}{c}
{\cal U}\\
{\cal A}\\
{\cal P}\\
{\cal V}\\
{\cal B}\\
{\cal Q}
\end{array}\right)\\
\ \\
\partial_2\left(\begin{array}{c}
{\cal U}\\
{\cal A}\\
{\cal P}\\
{\cal V}\\
{\cal B}\\
{\cal Q}
\end{array}\right)=
\left(\begin{array}{cccccc}
\frac{\partial_2p}{p} & 1 & 0 & 0 & 0 & 0\\
a & 0 & 1 & 0 & 0 & 0\\
0 & a & -\frac{\partial_2p}{p} & qb & 0 & -q\\
q & 0 & 0 & 0 & 0 & 0\\
0 & 0 & 0 & l & 0 & 0\\
-qb & 0 & 0 & 0 & l & 0
\end{array}\right)
\left(\begin{array}{c}
{\cal U}\\
{\cal A}\\
{\cal P}\\
{\cal V}\\
{\cal B}\\
{\cal Q}
\end{array}\right)
\end{array}
\label{lie1}
\end{equation}
which identically coincide with equations of the Lie sphere frame ---
see Appendix B. 
Let us define the wedge product $\wedge$ as in  Appendix C and
introduce the pseudo-Hermitian scalar product $(~, ~)$ in  $\Lambda^2{\bf
C}^4$ 
$$
({\bf a}\wedge {\bf b}, \ {\bf A}\wedge {\bf B})=
\frac{1}{2} 
 (({\bf a}, {\bf B})({\bf b}, {\bf A})-({\bf a}, {\bf A})({\bf b}, {\bf B}))
$$
(we hope that the same notation $(~,~)$ for pseudo-Hermitian scalar products
in
${\bf C}^4$ and $\Lambda^2{\bf C}^4$ will not lead to a confusion). A direct
computation shows that
the only nonzero products among the vectors
${\cal U}, {\cal A}, {\cal P}, {\cal V}, {\cal B}, {\cal Q}$ are
\begin{equation}
({\cal U}, \ {\cal P})=({\cal P}, \ {\cal U})=({\cal V}, \ {\cal Q})=
({\cal Q}, \ {\cal V})=-1, ~~~ ({\cal A}, \ {\cal A})=({\cal B}, \ {\cal
B})=1.
\label{scal1}
\end{equation}
This invariant pseudo-Hermitian scalar product corresponds to the quadratic
integral
$$
{\cal A}\bar {\cal A}+{\cal B}\bar {\cal B}-{\cal U}\bar {\cal P}-{\cal
P}\bar {\cal U}
-{\cal V}\bar {\cal Q}-{\cal Q}\bar {\cal V}
$$
of system (\ref{lie1}). Similarly, we can define the complex scalar product
$\{\ , \ \}$
in $\Lambda^2({\bf C}^4)$  (see Appendix C):
$$
\{{\bf a}\wedge {\bf b}, \ {\bf A}\wedge {\bf B}\}=
\frac{1}{2} 
 ({\bf a}\wedge {\bf b} \wedge {\bf A}\wedge {\bf B})
$$
 A direct computation shows that
the only nonzero products among the vectors
${\cal U}, {\cal A}, {\cal P}$, ${\cal V}, {\cal B}, {\cal Q}$ are
\begin{equation}
\{{\cal U}, \ {\cal P}\}=\{{\cal P}, \ {\cal U}\}=\{{\cal V}, \ {\cal Q}\}=
\{{\cal Q}, \ {\cal V}\}=-1, ~~~ \{{\cal A}, \ {\cal A}\}=
\{{\cal B}, \ {\cal B}\}=1.
\label{scal2}
\end{equation}
This invariant complex scalar product corresponds to the quadratic integral
$$
{\cal A}^2+{\cal B}^2-2{\cal U}{\cal P}-2{\cal V}{\cal Q}
$$
of system (\ref{lie1}). Notice the important
relation between the pseudo-Hermitian
scalar product $(~,~)$ and the complex scalar product $\{~,~\}$ in
$\Lambda^2({\bf C}^4)$:
\begin{equation}
(\xi, \ \psi)=
\{\xi, \ \bar \psi\}
\label{scal3}
\end{equation}
for any $\xi, \psi$ in $\Lambda^2({\bf C}^4)$ (see Appendix 3).
Let us finally introduce  in $\Lambda^2({\bf C}^4)$ the real vectors
$$
{\bf U}={\cal U}+\bar{\cal U}=i\ {{\mbox{\boldmath $\psi$}}}\wedge
\partial_1{{\mbox{\boldmath $\psi$}}}-i \
\overline{{{\mbox{\boldmath $\psi$}}}\wedge \partial_1{{\mbox{\boldmath
$\psi$}}}}, ~~~~
{\bf V}={\cal V}+\bar{\cal V}={{\mbox{\boldmath $\psi$}}}\wedge
\partial_2{{\mbox{\boldmath $\psi$}}}+
\overline{{{\mbox{\boldmath $\psi$}}}\wedge \partial_2{{\mbox{\boldmath
$\psi$}}}}.
$$

\begin{theorem}
Vectors ${\bf U}$ and ${\bf V}$ have zero norm:
$$
({\bf U}, \  {\bf U})= ({\bf V}, \  {\bf V})=0.
$$
Moreover, the triple
${\bf U}, \ \partial_2 {\bf U}, \  \partial_2^2{\bf U}$
is orthogonal to the triple
${\bf V}, \ \partial_1 {\bf V}, \  \partial_1^2{\bf V}$. Hence, ${\bf U}$
and 
${\bf V}$ are curvature spheres of a surface.

\end{theorem}

\centerline{Proof:}

\noindent It readily follows from (\ref{lie1}) that the triple
${\bf U}, \ \partial_2 {\bf U}, \  \partial_2^2{\bf U}$ is equivalent to
${\cal U}+\bar{\cal U},  \ {\cal V}+\bar {\cal V}, \  {\cal P}+\bar {\cal
P}$.
Similarly, the triple ${\bf V}, \ \partial_1 {\bf V}, \  \partial_1^2{\bf
V}$
is equivalent to
${\cal V}+\bar{\cal V},  \ {\cal B}+\bar {\cal B}, \  {\cal Q}+\bar {\cal
Q}$.
Conditions $({\bf U}, \  {\bf U})= ({\bf V}, \  {\bf V})=0$ and the
orthogonality of both triples follow by virtue
of (\ref{scal1}), (\ref{scal2}) and (\ref{scal3}). Let us show, for
instance, that 
$({\bf U}, \  {\bf U})=0$:
$$
\begin{array}{c}
({\bf U}, \  {\bf U})=({\cal U} +\bar {\cal U}, \ {\cal U} +\bar {\cal U})=
({\cal U}, \ {\cal U})+({\cal U}, \ \bar {\cal U})+(\bar {\cal U}, \ {\cal
U})+
(\bar {\cal U}, \ \bar {\cal U})= \\
\ \\
({\cal U}, \ {\cal U})+\{{\cal U}, \ {\cal U}\}+
\overline {({\cal U}, \ \bar {\cal U})}+
\{\bar {\cal U}, \  {\cal U}\}=
({\cal U}, \ {\cal U})+\{{\cal U}, \ {\cal U}\}+
\overline {\{{\cal U}, \  {\cal U}\}}+
\{ {\cal U}, \  \bar {\cal U}\}= \\
\ \\
({\cal U}, \ {\cal U})+\{{\cal U}, \ {\cal U}\}+
\overline {\{{\cal U}, \  {\cal U}\}}+
( {\cal U}, \   {\cal U}),
\end{array}
$$
which is zero by virtue of (\ref{scal1}) and  (\ref{scal2}).

\medskip

\noindent {\bf Remark 1.} It is straightforward to show that any  surface
can 
be obtained (locally) by a construction described above.

\medskip

\noindent {\bf Remark 2.} In view of  (\ref{norm}) the surface
${\mbox{\boldmath $\psi$}}(R^1, R^2) \in {\bf CP}^3$ defines a Legendre
submanifold of the quadric $({\mbox{\boldmath $\psi$}}, \
{\mbox{\boldmath $\psi$}})=0$ equipped with a real contact form
$i({\mbox{\boldmath $\psi$}}, \ d{\mbox{\boldmath $\psi$}})$. I would like
to thank L. Mason 
for clarifying this point.

\section{Surfaces possessing 3-parameter families of Lie deformations and
commuting 
Schr\"odinger operators with magnetic fields}

In contrast with the Euclidean geometry, where a surface is uniquely
determined by its 
first and second fundamental forms, there exist examples of surfaces in Lie
sphere
geometry which are not uniquely specified by the Lie-invariant metric
$$
-pq \ dR^1dR^2
$$
and the conformal class of the cubic form
$$
p \, (dR^1)^3-q \, (dR^2)^3.
$$
Such surfaces are called Lie-applicable (Lie-deformable).
In this section we consider examples of surfaces possessing 3-parameter
families of Lie 
deformations. A calculation
similar to the one done by Finikov in \cite{Finikov37} shows that these
surfaces 
are characterized by the constraints
\begin{equation}
\partial_1\partial_2\ln p = c\ pq, ~~~ \partial_1\partial_2\ln q = c\ pq,
\label{c}
\end{equation}
where $c$ is a constant. There are different
cases to distinguish depending on the value of $c$. Here we discuss the two
simplest 
cases $c=0$ and $c=1$ (for $c\ne 0, 1$ the formulae become more
complicated).

\medskip

\noindent {\bf Case  c=0} Utilizing transformations (\ref{newR}) and
(\ref{newpq}),
 we can represent $p$ and $q$ in the form
$$
p=\psi_1'(R^1), ~~~ q=-\psi_2'(R^2),
$$
implying, after the substitution into (\ref{VW}) and  elementary
integration, the following
expressions for $V$ and $ W$
$$
V=\epsilon_1+\epsilon_0\psi_1-\psi_2\psi_1''-\frac{1}{2}\rho_2\psi_1^2, ~~~
W=\epsilon_2+\epsilon_0\psi_2-\psi_1\psi_2''-\frac{1}{2}\rho_1\psi_2^2.
$$
Moreover,  $\psi_1$ and $ \psi_2$ satisfy the ODE's
$$
\psi_1''=\alpha \psi_1^2 +\rho_1 \psi_1+s_1, ~~~
\psi_2''=\alpha \psi_2^2 +\rho_2 \psi_2+s_2
$$
implying that $\psi_1$ and $\psi_2$ are elliptic functions.
Here $\epsilon_0, \epsilon_1, \epsilon_2, \alpha, \rho_1, \rho_2, s_1, s_2$
are arbitrary constants
(if $\alpha$ is nonzero one can always reduce $\rho_1$ and $\rho_2$ to zero
by adding constants to $\psi_1, \psi_2$). Notice that for given $p$ and $q$
the corresponding $V$ and $W$ are
determined up to three arbitrary constants $\epsilon_0, \epsilon_1,
\epsilon_2$, which
are thus responsible for Lie deformations. It is important to emphasize
the linear dependence of
$V$ and $W$ on the deformation parameters. This readily follows from
 (\ref{VW}), indeed, for given $p$ and $q$ these equations are linear in $V$
and $W$.
 The corresponding system  (\ref{LieWil}) takes the form
\begin{equation}
\begin{array}{c}
\partial_1^2{{\mbox{\boldmath $\psi$}}}=
-i \psi_1' \ \partial_2{{\mbox{\boldmath $\psi$}}}+\frac{1}{2}
(\epsilon_1+\epsilon_0\psi_1-\psi_2\psi_1''-\frac{1}{2}\rho_2\psi_1^2) \
{{\mbox{\boldmath $\psi$}}}, \\
\ \\
\partial_2^2{{\mbox{\boldmath $\psi$}}}=
-i \psi_2' \ \partial_1{{\mbox{\boldmath $\psi$}}}+\frac{1}{2}
(\epsilon_2+\epsilon_0\psi_2-\psi_1\psi_2''-
\frac{1}{2}\rho_1\psi_2^2) \ {{\mbox{\boldmath $\psi$}}},
\end{array}
\label{c=0}
\end{equation}
which, upon the addition and subtraction,  readily rewrites
in the form
$$
H{{\mbox{\boldmath $\psi$}}}=\lambda \ {{\mbox{\boldmath $\psi$}}}, ~~~
F{{\mbox{\boldmath $\psi$}}}=\mu \ {{\mbox{\boldmath $\psi$}}}.
$$
Here $H$ and $F$ are commuting Schr\"odinger operators with magnetic terms
$$
H=\left(\partial_1+\frac{i}{2}\psi_2'\right)^2+
\left(\partial_2+\frac{i}{2}\psi_1'\right)^2+V_H,
$$
$$
F=\left(\partial_1-\frac{i}{2}\psi_2'\right)^2-
\left(\partial_2-\frac{i}{2}\psi_1'\right)^2+V_F,
$$
$\lambda = \frac{\epsilon_1+\epsilon_2}{2}$ and
$\mu = \frac{\epsilon_1-\epsilon_2}{2}$ are the eigenvalues, and the
potentials $V_H$ and $V_F$ are given by
$$
V_H=\frac{1}{4}(2\psi_2\psi_1''+2\psi_1\psi_2''+\rho_2\psi_1^2+\rho_1\psi_2^
2
+(\psi_2')^2+(\psi_1')^2-2\epsilon_0(\psi_1+\psi_2)),
$$
$$
V_F=\frac{1}{4}(2\psi_2\psi_1''-2\psi_1\psi_2''+\rho_2\psi_1^2-\rho_1\psi_2^
2
+(\psi_2')^2-(\psi_1')^2-2\epsilon_0(\psi_1-\psi_2)).
$$
For generic values of constants, operators $H$ and $F$ will be non-singular
and 
doubly periodic.
The spectral theory of these operators will be discussed elsewhere.
\bigskip

\noindent {\bf Case  c=1} Here
$$
\partial_1\partial_2\ln p =  pq, ~~~ \partial_1\partial_2\ln q =  pq,
$$
implying
$$
p=\frac{1}{R^2-R^1}\sqrt {\frac{f_2}{f_1}}, ~~~
q=\frac{1}{R^1-R^2}\sqrt {\frac{f_1}{f_2}},
$$
where $f_1$ and $f_2$ are functions of $R^1$ and $R^2$, respectively.
The corresponding $V$ and $W$ are given by
$$
\begin{array}{c}
V=\partial_1^2(\ln q) +\frac{1}{2}(\partial_1 q/q)^2-
\frac{\epsilon_0+\epsilon_1R^1+\epsilon_2(R^1)^2}{f_1}, \\
\  \\
W=\partial_2^2(\ln p) +\frac{1}{2}(\partial_2 p/p)^2+
\frac{\epsilon_0+\epsilon_1R^2+\epsilon_2(R^2)^2}{f_2},
\end{array}
$$
where the constants $\epsilon_0, \epsilon_1, \epsilon_2$ are responsible for
 Lie deformations. These surfaces are known to
 have both families of curvature lines in linear complexes \cite{Blaschke}.
Introducing the rescaled vector ${\bf R}$ by the formula
$$
{{\mbox{\boldmath $\psi$}}}=(f_1f_2)^{1/4}{\bf R},
$$
we readily rewrite  equations  (\ref{LieWil}) in the
equivalent form
\begin{equation}
\begin{array}{c}
f_1\ \partial_1^2{\bf R}+\frac{1}{2}f_1'\ \partial_1{\bf R}+
i\frac{\sqrt {f_1f_2}}{R^2-R^1}\ \partial_2{\bf R}= \\
\ \\
\left( \frac{3f_1}{4(R^2-R^1)^2}+\frac{f_1'}{4(R^2-R^1)}-
\frac{i}{2}\frac{\sqrt {f_1f_2}}{(R^2-R^1)^2}-
\frac{\epsilon_0+\epsilon_1 R^1+\epsilon_2 ({R^1})^2}{2} \right){\bf R},\\
\ \\
f_2\ \partial_2^2{\bf R}+\frac{1}{2}f_2'\ \partial_2{\bf R}+
i\frac{\sqrt {f_1f_2}}{R^2-R^1}\ \partial_1{\bf R}=\\
\ \\
\left( \frac{3f_2}{4(R^2-R^1)^2}-\frac{f_2'}{4(R^2-R^1)}+
\frac{i}{2}\frac{\sqrt {f_1f_2}}{(R^2-R^1)^2}+
\frac{\epsilon_0+\epsilon_1 R^2+\epsilon_2 ({R^2})^2}{2} \right){\bf R},\\
\end{array}
\label{bfR}
\end{equation}
where $f_1'$ and $f_2'$ denote the derivatives of
$f_1(R^1)$ and $f_2(R^2)$, respectively.
Solving for $\frac{\epsilon_1}{2}{\bf R}$ and $\frac{\epsilon_0}{2}{\bf R}$,
we
arrive at the eigenfunction equations
$$
 H {\bf R}+\frac{\epsilon_1}{2}{\bf R}=0, ~~~~
 F {\bf R}+\frac{\epsilon_0}{2}{\bf R}=0
$$
where the operator $ H$ is of the form
\begin{equation}
\sqrt {g^{11}g^{22}}\left(i\partial_x +A\right)
\frac{g^{11}}{\sqrt{g^{11}g^{22}}}
\left(i\partial_x +A\right)+
\sqrt {g^{11}g^{22}}\left(i\partial_y +B\right)
\frac{g^{22}}{\sqrt{g^{11}g^{22}}}
\left(i\partial_y +B\right)+h.
\label{Laplace}
\end{equation}
Here $g^{11}$ and $g^{22}$ are the components of a diagonal metric of
St\"ackel
type
$$
g^{11}=\frac{f_1}{R^2-R^1}, ~~~ g^{22}=\frac{f_2}{R^2-R^1},
$$
$A$ and $B$ are the components of the magnetic vector potential
$$
A=-\frac{1}{2(R^2-R^1)}\sqrt{\frac{f_2}{f_1}}, ~~~
B=-\frac{1}{2(R^2-R^1)}\sqrt{\frac{f_1}{f_2}},
$$
and $h$ is the scalar potential
$$
h=\frac{f_1'-f_2'}{4(R^2-R^1)^2}+\frac{f_1+f_2}{2(R^2-R^1)^3}+
\frac{\epsilon_2}{2}(R^1+R^2).
$$
Geometrically, $ H$ represents the Laplace-Beltrami operator corresponding
to
the metric $g^{11}, g^{22}$ in the
magnetic potential $A, B$ and the scalar potential $h$.
Notice that the scalar potential $h$ can be represented  in a simple
coordinate-free form
$$
h=K+\frac{\epsilon_2}{2}(R^1+R^2)
$$
where $K$ is the Gaussian curvature of the metric $g^{11}, \ g^{22}$. The
second term $R^1+R^2$ is nothing but the trace of the Killing tensor of the
St\"ackel metric $g^{11}, \ g^{22}$, and hence also makes an invariant
sense.
Computation of the magnetic field implies
$$
(\partial_1 B-\partial_2 A)\ dR^1\wedge dR^2=
-\left(\frac{f_1'-f_2'}{4(R^2-R^1)^2}+\frac{f_1+f_2}{2(R^2-R^1)^3} \right)\
 \frac{R^2-R^1}{\sqrt {f_1f_2}} \ dR^1\wedge dR^2=
-K\ d\sigma
$$
where
$$
d\sigma=\frac{R^2-R^1}{\sqrt {f_1f_2}} \ dR^1\wedge dR^2
$$
is the area form of the metric $g^{11}, \ g^{22}$.
Thus, the magnetic field also makes an invariant sense. Notice that in the
case when
$$
f^1=4(R^1)^3+a(R^1)^2+bR^1+c, ~~~ f^2=-4(R^2)^3-a(R^2)^2-bR^2-c
$$
are cubic polynomials, the Gaussian curvature $K=1$ and the operator $H$
represents Dirac monopole on the unit sphere in the spherical-conical
coordinates $R^1, R^2$.
The scalar potential $h$ (which in this case is proportional to $R^1+R^2$)
has a
 meaning of the external quadratic potential.
We refer to \cite{Fer9} for the discussion of some algebraic
aspects of spectral theory  of such operators in the particular
 case $\epsilon_2=0$.
The general situation will be discussed elsewhere.

The examples discussed in this section clearly demonstrate that there exists
a one-to-one
correspondence between commuting Schr\"odinger operators with magnetic
fields and
Lie-applicable surfaces which possess multi-parameter families of Lie
deformations.

\section{Canal surfaces}

Our approach to the canal surfaces will be based on a linear system
\begin{equation}
\begin{array}{c}
\partial_1^2{{\mbox{\boldmath $\psi$}}}=
-i \ p \ \partial_2{{\mbox{\boldmath $\psi$}}}+
\frac{1}{2}(V+i \ \partial_2p) \ {{\mbox{\boldmath $\psi$}}} \\
\ \\
\partial_2^2{{\mbox{\boldmath $\psi$}}}=\frac{1}{2}W \
{{\mbox{\boldmath $\psi$}}} \\
\end{array}
\label{canal}
\end{equation}
which is a specialization of (\ref{LieWil}) corresponding to $q=0$.
The compatibility conditions of system (\ref{canal}) take the form
\begin{equation}
\begin{array}{c}
\partial_2^3p-2W\partial_2p-p\partial_2 W=0 \\
\ \\
\partial_1W=\partial_2V=0.
\end{array}
\label{canal1}
\end{equation}
Since $q=0$, we cannot use formulae  (\ref{framelie}). Instead, we
 introduce the  four vectors
\begin{equation}
{{\mbox{\boldmath $\psi$}}}, ~~~
{{\mbox{\boldmath $\psi$}}}_1=\partial_1{{\mbox{\boldmath $\psi$}}}, ~~~
{{\mbox{\boldmath $\psi$}}}_2=\partial_2{{\mbox{\boldmath
$\psi$}}}-\frac{1}{2}\frac{\partial_2p}{p}{{\mbox{\boldmath $\psi$}}}, ~~~
{\mbox{\boldmath $\eta$}}=\partial_1\partial_2{{\mbox{\boldmath $\psi$}}}-
\frac{1}{2}\frac{\partial_2p}{p}\partial_1{{\mbox{\boldmath $\psi$}}}
\label{framelie1}
\end{equation}
which satisfy the linear system
\begin{equation}
\begin{array}{c}
\partial_1\left(\begin{array}{c}
{{\mbox{\boldmath $\psi$}}}\\
{{\mbox{\boldmath $\psi$}}}_1\\
{{\mbox{\boldmath $\psi$}}}_2\\
{\mbox{\boldmath $\eta$}}
\end{array}\right)=
\left(\begin{array}{cccc}
0& 1&0&0\\
\frac{1}{2}b & 0 & -ip &0\\
\frac{1}{2}k&0&0&1\\
-\frac{i}{2}p
a&\frac{1}{2}k&\frac{1}{2}b&0
\end{array}\right)
\left(\begin{array}{c}
{{\mbox{\boldmath $\psi$}}}\\
{{\mbox{\boldmath $\psi$}}}_1\\
{{\mbox{\boldmath $\psi$}}}_2\\
{\mbox{\boldmath $\eta$}}
\end{array}\right)\\
\ \\
\partial_2\left(\begin{array}{c}
{{\mbox{\boldmath $\psi$}}}\\
{{\mbox{\boldmath $\psi$}}}_1\\
{{\mbox{\boldmath $\psi$}}}_2\\
{\mbox{\boldmath $\eta$}}
\end{array}\right)=
\left(\begin{array}{cccc}
\frac{1}{2}\frac{\partial_2p}{p}& 0&1&0\\
0 & \frac{1}{2}\frac{\partial_2p}{p}& 0&1\\
\frac{1}{2}a&0 &-\frac{1}{2}\frac{\partial_2p}{p}&0\\
0&\frac{1}{2}a
&0&-\frac{1}{2}\frac{\partial_2p}{p}
\end{array}\right)\left(
\begin{array}{c}
{{\mbox{\boldmath $\psi$}}}\\
{{\mbox{\boldmath $\psi$}}}_1\\
{{\mbox{\boldmath $\psi$}}}_2\\
{\mbox{\boldmath $\eta$}}
\end{array}\right)
\end{array}
\label{lieframe1}
\end{equation}
where  the notation
$$
k=- \partial_1\partial_2(\ln p), ~~~
a=W-\partial_2^2(\ln p)-\frac{1}{2}(\partial_2 \ln p)^2, ~~~
b=V
$$
is introduced.
Compatibility conditions of equations (\ref{lieframe1}) imply
$$
\partial_1\partial_2 \ln p= -k, ~~~
\partial_1a=\partial_2k+\frac{\partial_2p}{p}\, k, ~~~
\partial_2b=0, ~~~
p\,  \partial_2a+2\, a\, \partial_2p  = 0.
$$
An important property of  system (\ref{lieframe1}) is the existence of the
quadratic integral
\begin{equation}
-{{\mbox{\boldmath $\psi$}}}\bar {\mbox{\boldmath $\eta$}}+{{\mbox{\boldmath
$\psi$}}}_1{\bar {{\mbox{\boldmath $\psi$}}}}_2+
{{\mbox{\boldmath $\psi$}}}_2 {\bar {{\mbox{\boldmath $\psi$}}}}_1 -
{\mbox{\boldmath $\eta$}}{\bar {{\mbox{\boldmath $\psi$}}}}
\label{integral11}
\end{equation}
which defines an invariant pseudo-Hermitian scalar product of the signature
$(2, 2)$ on the space of solutions of system (\ref{lieframe1}).
Using  (\ref{framelie1}), this integral can be rewritten in the form
\begin{equation}
-{{\mbox{\boldmath $\psi$}}} \partial_1\partial_2{\bar{{\mbox{\boldmath
$\psi$}}}}+ 
\partial_1{{\mbox{\boldmath $\psi$}}} \partial_2 {\bar {{\mbox{\boldmath
$\psi$}}}}+
\partial_2{{\mbox{\boldmath $\psi$}}} \partial_1{\bar {{\mbox{\boldmath
$\psi$}}}}-
{\bar {{\mbox{\boldmath $\psi$}}}} \partial_1\partial_2{{{\mbox{\boldmath
$\psi$}}}}.
\label{integral21}
\end{equation}
The  invariant pseudo-Hermitian scalar product (\ref{integral11})
implies the existence of a basis of solutions of (\ref{lieframe1})
such that 
\begin{equation}
({{\mbox{\boldmath $\psi$}}}_1, \ {{\mbox{\boldmath
$\psi$}}}_2)=({{\mbox{\boldmath $\psi$}}}_2, \ {{\mbox{\boldmath
$\psi$}}}_1)=1, ~~~~
({{\mbox{\boldmath $\psi$}}}, \ {\mbox{\boldmath $\eta$}})=({\mbox{\boldmath
$\eta$}}, 
\ {{\mbox{\boldmath $\psi$}}})=-1,
\label{Herm1}
\end{equation}
all other scalar products being zero.
Equations (\ref{Herm1}) are obviously equivalent to
$$
(\partial_1 {{\mbox{\boldmath $\psi$}}}, \partial_2 {{\mbox{\boldmath
$\psi$}}})=
(\partial_2 {{\mbox{\boldmath $\psi$}}}, \partial_1 {{\mbox{\boldmath
$\psi$}}})=1, ~~~
(\partial_1 \partial_2 {{\mbox{\boldmath $\psi$}}}, {{\mbox{\boldmath
$\psi$}}})=
({{\mbox{\boldmath $\psi$}}}, \partial_1 \partial_2 {{\mbox{\boldmath
$\psi$}}})=-1.
$$
Here $(~,~)$ denotes the pseudo-Hermitian scalar product in
${\bf C}^4$ of the signature $(2, \ 2)$ as in  Appendix C.
Equations (\ref{lieframe1}) also imply that
the determinant
${{\mbox{\boldmath $\psi$}}}\wedge {{\mbox{\boldmath $\psi$}}}_1 \wedge
{{\mbox{\boldmath $\psi$}}}_2\wedge {\mbox{\boldmath $\eta$}}$
is invariant:
$$
\partial_1 ({{\mbox{\boldmath $\psi$}}}\wedge {{\mbox{\boldmath $\psi$}}}_1
\wedge {{\mbox{\boldmath $\psi$}}}_2\wedge {\mbox{\boldmath $\eta$}})=
\partial_2 ({{\mbox{\boldmath $\psi$}}}\wedge {{\mbox{\boldmath $\psi$}}}_1
 \wedge {{\mbox{\boldmath $\psi$}}}_2\wedge {\mbox{\boldmath $\eta$}})=0.
$$
 Thus, besides 
(\ref{Herm1}), we can impose the additional constraint
\begin{equation}
{{\mbox{\boldmath $\psi$}}}\wedge {{\mbox{\boldmath $\psi$}}}_1 \wedge
{{\mbox{\boldmath $\psi$}}}_2\wedge {\mbox{\boldmath $\eta$}}=1.
\label{det1}
\end{equation}
From now on, we fix a null-tetrad
${{\mbox{\boldmath $\psi$}}}, \ {{\mbox{\boldmath $\psi$}}}_1, \
{{\mbox{\boldmath $\psi$}}}_2, \ {\mbox{\boldmath $\eta$}}$
satisfying both (\ref{Herm1}) and (\ref{det1}). Notice that such a basis
 is defined up to the action of the group
$SU(2, 2)$ which preserves both (\ref{Herm1}) and (\ref{det1}).

Introducing the basis in $\Lambda^2{\bf C}^4$ as follows
$$
\begin{array}{c}
{\cal U}=i \ {{\mbox{\boldmath $\psi$}}} \wedge {{\mbox{\boldmath
$\psi$}}}_1, ~~~ {\cal V}={{\mbox{\boldmath $\psi$}}} \wedge
{{\mbox{\boldmath $\psi$}}}_2, \\
\ \\
{\cal A} =i \ {{\mbox{\boldmath $\psi$}}}_2\wedge {{\mbox{\boldmath
$\psi$}}}_1+i \ {{\mbox{\boldmath $\psi$}}}\wedge {\mbox{\boldmath $\eta$}},
~~~
{\cal B} ={{\mbox{\boldmath $\psi$}}}_1\wedge {{\mbox{\boldmath
$\psi$}}}_2+{{\mbox{\boldmath $\psi$}}}\wedge {\mbox{\boldmath $\eta$}}, \\
\ \\
{\cal P} = 2i\, {{\mbox{\boldmath $\psi$}}}_2\wedge {\mbox{\boldmath
$\eta$}}, ~~~
{\cal Q} = 2\, {{\mbox{\boldmath $\psi$}}}_1\wedge {\mbox{\boldmath
$\eta$}},
\end{array}
$$
we arrive at the  equations
\begin{equation}
\begin{array}{c}
\partial_1\left(\begin{array}{c}
{\cal U}\\
{\cal A}\\
{\cal P}\\
{\cal V}\\
{\cal B}\\
{\cal Q}
\end{array}\right)=
\left(\begin{array}{cccccc}
0 & 0 & 0 & p & 0 & 0\\
k & 0 & 0 & 0 & 0 & 0\\
0 & k & 0 & -pa & 0 & 0\\
0 & 0 & 0 & 0 & 1 & 0\\
0 & 0 & 0 & b & 0 & 1\\
pa & 0 & -p & 0 & b & 0
\end{array}\right)
\left(\begin{array}{c}
{\cal U}\\
{\cal A}\\
{\cal P}\\
{\cal V}\\
{\cal B}\\
{\cal Q}
\end{array}\right)\\
\ \\
\partial_2\left(\begin{array}{c}
{\cal U}\\
{\cal A}\\
{\cal P}\\
{\cal V}\\
{\cal B}\\
{\cal Q}
\end{array}\right)=
\left(\begin{array}{cccccc}
\frac{\partial_2p}{p} & 1 & 0 & 0 & 0 & 0\\
a & 0 & 1 & 0 & 0 & 0\\
0 & a & -\frac{\partial_2p}{p} & 0 & 0 & 0\\
0 & 0 & 0 & 0 & 0 & 0\\
0 & 0 & 0 & 0 & 0 & 0\\
0 & 0 & 0 & 0 & 0 & 0
\end{array}\right)
\left(\begin{array}{c}
{\cal U}\\
{\cal A}\\
{\cal P}\\
{\cal V}\\
{\cal B}\\
{\cal Q}
\end{array}\right).
\end{array}
\label{lie11}
\end{equation}
Let us define the pseudo-Hermitian scalar product $(~, ~)$ and
the complex scalar product
$\{\ , \ \}$ in 
 $\Lambda^2{\bf C}^4$ as in Appendix 3.
 A direct computation shows that
the only nonzero products among the vectors
${\cal U}, {\cal A}, {\cal P}, {\cal V}, {\cal B}, {\cal Q}$ are
\begin{equation}
({\cal U}, \ {\cal P})=({\cal P}, \ {\cal U})=({\cal V}, \ {\cal Q})=
({\cal Q}, \ {\cal V})=-1, ~~~ ({\cal A}, \ {\cal A})=({\cal B}, \ {\cal
B})=1
\label{scal11}
\end{equation}
and
\begin{equation}
\{{\cal U}, \ {\cal P}\}=\{{\cal P}, \ {\cal U}\}=\{{\cal V}, \ {\cal Q}\}=
\{{\cal Q}, \ {\cal V}\}=-1, ~~~ \{{\cal A}, \ {\cal A}\}=
\{{\cal B}, \ {\cal B}\}=1,
\label{scal21}
\end{equation}
respectively. Thus,
\begin{equation}
(\xi, \ \psi)=
\{\xi, \ \bar \psi\}
\label{scal31}
\end{equation}
for any $\xi, \psi$ in $\Lambda^2({\bf C}^4)$.
Introducing  in $\Lambda^2({\bf C}^4)$ the real vectors
$$
{\bf U}={\cal U}+\bar{\cal U}=i\ {{\mbox{\boldmath $\psi$}}}\wedge
\partial_1{{\mbox{\boldmath $\psi$}}}-i \
\overline{{{\mbox{\boldmath $\psi$}}}\wedge \partial_1{{\mbox{\boldmath
$\psi$}}}}, ~~~~
{\bf V}={\cal V}+\bar{\cal V}={{\mbox{\boldmath $\psi$}}}\wedge
\partial_2{{\mbox{\boldmath $\psi$}}}+
\overline{{{\mbox{\boldmath $\psi$}}}\wedge
\partial_2{{\mbox{\boldmath $\psi$}}}},
$$
we can formulate the main result of this section

\begin{theorem}
Vectors ${\bf U}$ and ${\bf V}$ have zero norm:
$$
({\bf U}, \  {\bf U})= ({\bf V}, \  {\bf V})=0.
$$
Moreover, the triple
${\bf U}, \ \partial_2 {\bf U}, \  \partial_2^2{\bf U}$
is orthogonal to the triple
${\bf V}, \ \partial_1 {\bf V}, \  \partial_1^2{\bf V}$. Hence, ${\bf U}$
and 
${\bf V}$ are curvature spheres of a surface.
This surface will be a canal surface since
$\partial_2{\bf V}=0$.
Any canal surface can be obtained (locally) by this construction.

\end{theorem}

\medskip

The proof of this theorem copies the proof of Theorem 1 from section 2:
it readily follows from (\ref{lie11}) that the triple
${\bf U}, \ \partial_2 {\bf U}, \  \partial_2^2{\bf U}$ is equivalent to
${\cal U}+\bar{\cal U},  \ {\cal V}+\bar {\cal V}, \  {\cal P}+\bar {\cal
P}$.
Similarly, the triple ${\bf V}, \ \partial_1 {\bf V}, \  \partial_1^2{\bf
V}$
is equivalent to
${\cal V}+\bar{\cal V},  \ {\cal B}+\bar {\cal B}, \  {\cal Q}+\bar {\cal
Q}$.
The conditions $({\bf U}, \  {\bf U})= ({\bf V}, \  {\bf V})=0$
and the orthogonality of both triples follow by virtue
of (\ref{scal11}), (\ref{scal21}) and (\ref{scal31}).

\bigskip

{\bf Example.} Let us consider the Landau operator
$$
H=\frac{1}{2}(i \ \partial_x -My)^2+\frac{1}{2}(i \ \partial_y)^2
$$
describing a quantum particle in the homogeneous magnetic field
($M={\rm const}$). 
Operator $H$ obviously commutes with the operator
$$
F=-\partial_x^2
$$
so that the equations for their joint eigenfunctions
$$
H {\mbox{\boldmath $\psi$}}=\lambda \ {\mbox{\boldmath $\psi$}}, ~~~~
F{\mbox{\boldmath $\psi$}} = k^2 \ {\mbox{\boldmath $\psi$}}
$$
can be rewritten in the form
\begin{equation}
\begin{array}{c}
{\mbox{\boldmath $\psi$}}_{yy}=-2iMy \ {\mbox{\boldmath
$\psi$}}_x+(M^2y^2+k^2-2\lambda) \ {\mbox{\boldmath $\psi$}} \\
\ \\
{\mbox{\boldmath $\psi$}}_{xx}=-k^2 \ {\mbox{\boldmath $\psi$}}.
\label{Landau}
\end{array}
\end{equation}
System (\ref{Landau}) is obviously of the form (\ref{canal})
under the identification $p=2MR^1, \ V=2M^2(R^1)^2 +2k^2-4\lambda, \
W=-2k^2, \ y=R^1, \ x=R^2$. In what follows we assume $M=1$ (this can be
achieved by a rescaling $\tilde x=x\sqrt M, \ \tilde y=y\sqrt M, \
 \tilde k=k/\sqrt M, \ \tilde \lambda=\lambda / M$). The corresponding
linear system for
${\mbox{\boldmath $\psi$}}$ possesses 4 linearly independent solutions
$$
e^{ikx}\psi_1(y+k), ~~~ e^{-ikx}\psi_1(y-k),
~~~ e^{ikx}\psi_2(y+k), ~~~ e^{-ikx}\psi_2(y-k)
$$
where $\psi_1, \ \psi _2$ form a basis of solutions of Hermite's equation
$$
\psi ''=(y^2-2\lambda)\psi.
$$
Let us introduce the complex 4-vector
$$
{{\mbox{\boldmath $\psi$}}}=\left(e^{ikx}\psi_1(y+k), \ e^{-ikx}\psi_1(y-k),
\
\frac{e^{-ikx}\psi_2(y-k)}{2ikW}, \ \frac{e^{ikx}\psi_2(y+k)}{2ikW}\right)
$$
where $W=\psi _1\psi _2'-\psi _2\psi _1'={\rm const}$ is the Wronskian.
Then
$$
({{\mbox{\boldmath $\psi$}}}_x, \ {{\mbox{\boldmath
$\psi$}}}_y)=({{\mbox{\boldmath $\psi$}}}_y, \ {{\mbox{\boldmath
$\psi$}}}_x)=1, ~~~
({{\mbox{\boldmath $\psi$}}}, \ {{\mbox{\boldmath $\psi$}}}_{xy})=
({{\mbox{\boldmath $\psi$}}}_{xy}, \ {{\mbox{\boldmath $\psi$}}})=-1,
$$
all other products being zero. Moreover,
$$
{{\mbox{\boldmath $\psi$}}} \wedge {{\mbox{\boldmath $\psi$}}}_x \wedge
{{\mbox{\boldmath $\psi$}}}_y \wedge {{\mbox{\boldmath $\psi$}}}_{xy}= -1.
$$ 
A direct computation implies
$$
{\bf V}={{\mbox{\boldmath $\psi$}}}\wedge {{\mbox{\boldmath $\psi$}}}_x+
\overline{{{\mbox{\boldmath $\psi$}}}\wedge {{\mbox{\boldmath $\psi$}}}_x}=
(y^0, \ y^1, \ y^2, \ y^3, \ y^4, \ y^5)
$$
where
$$
y^0=\frac{\psi_1(y-k)\psi_2(y+k)-\psi_2(y-k)\psi_1(y+k)}{W},
$$
$$
y^1=-\frac{\psi_1(y-k)\psi_2(y+k)+\psi_2(y-k)\psi_1(y+k)}{W},
$$
$$
y^2=y^3=0,
$$
$$
y^4=-2k\psi_1(y-k)\psi_1(y+k)+\frac{\psi_2(y-k)\psi_2(y+k)}{2kW^2},
$$
$$
y^5=-2k\psi_1(y-k)\psi_1(y+k)-\frac{\psi_2(y-k)\psi_2(y+k)}{2kW^2}.
$$
Obviously, $-(y^0)^2+(y^1)^2+(y^2)^2+(y^3)^2+(y^4)^2-(y^5)^2=0$.
Dividing by $y^0+y^1$, we obtain the normalised vector
$$
{\bf V}={{\mbox{\boldmath $\psi$}}}\wedge {{\mbox{\boldmath $\psi$}}}_x+
\overline{{{\mbox{\boldmath $\psi$}}}\wedge {{\mbox{\boldmath $\psi$}}}_x}=
\left(\frac{y^0}{y^0+y^1}, \ \frac{y^1}{y^0+y^1}, \ 0, \ 0, \
\frac{y^4}{y^0+y^1}, \ \frac{y^5}{y^0+y^1}\right).
$$
Since the centers of the corresponding spheres lie on the $z$-axis, our
surface 
is a surface of revolution.
Parametric equations of  centers $z$ and radii $R$ are
$$
z=\frac{y^4}{y^0+y^1}=kW\frac{\psi_1(y-k)}{\psi_2(y-k)}-
\frac{1}{4kW} \frac{\psi_2(y+k)}{\psi_1(y+k)},
$$
$$
R=\frac{y^5}{y^0+y^1}=kW\frac{\psi_1(y-k)}{\psi_2(y-k)}+
\frac{1}{4kW} \frac{\psi_2(y+k)}{\psi_1(y+k)}.
$$
The Landau levels correspond to $\lambda =\frac{2n+1}{2}$. For $n=0$ we have
$$
\psi_1=e^{-\frac{y^2}{2}}, ~~~
\psi_2=e^{-\frac{y^2}{2}}\int \limits _0^y e^{\xi^2}\ d \xi, ~~~ W=1.
$$
I would like to thank
K.R. Khusnutdinova for the investigation of this and other examples.
The details will be given elsewhere.

\section{Appendix A: Wilczynski's projective frame}

Based on \cite{Wilczynski} (see also \cite{Bol}, \cite{Lane},
\cite{Fer8}, \cite{Fer002},
\cite{Sasaki1}), let us briefly recall
the standard way of defining surfaces $M^2$
in projective space $P^3$ in terms of solutions of a linear system
(\ref{r}) satisfying the compatibility conditions (\ref{GC1}). For any fixed
$\beta, \gamma, V, W$ satisfying (\ref{GC1}), the linear system (\ref{r})
is compatible and possesses a solution ${\bf r}=(r^0, r^1, r^2, r^3)$
where $r^i(x, y)$ can be regarded as homogeneous coordinates of a surface in
projective space $P^3$. In what follows, we assume that our surfaces are
hyperbolic 
and the corresponding asymptotic coordinates $x$ and $y$ are real.
Even though the coefficients $\beta, \gamma, V, W$ define a surface $M^2$
uniquely up to projective equivalence via (\ref{r}),
it is not entirely correct
to regard $\beta, \gamma, V, W$ as projective invariants. Indeed, the
asymptotic
coordinates $x, y$ are only defined up to an arbitrary reparametrization
of the form
\begin{equation}
x^*=f(x), ~~~~ y^*=g(y)
\label{newxy}
\end{equation}
which induces a scaling of the surface vector according to
\begin{equation}
{\bf r}^*=\sqrt {f'(x)g'(y)}~{\bf r}.
\label{newr}
\end{equation}
Thus \cite[p.\ 1]{Bol}, the form of
equations (\ref{r}) is preserved by the above transformation
with the new coefficients $\beta^*, \gamma^*, V^*, W^*$ given by
\begin{equation}
\begin{array}{c}
\beta^{*}=\beta g'/(f')^2, ~~~~  V^{*}(f')^2=V+S(f)\\
\ \\
\gamma^{*}=\gamma f'/(g')^2, ~~~~ W^{*}(g')^2=W+S(g),
\end{array}
\label{newbetagamma}
\end{equation}
where $S(\,\cdot\,)$ is the  Schwarzian derivative, that is
$$
S(f)=\frac{f'''}{f'} - \frac{3}{2} \left(\frac{f''}{f'}\right)^2.
$$
The transformation formulae (\ref{newbetagamma}) imply that the symmetric
2-form
$$
2 \beta \gamma\,dxdy
$$
and the conformal class of the cubic form
$$
\beta \,dx^3+\gamma \,dy^3
$$
are absolute projective invariants. They are known as the projective metric
and
the Darboux cubic form, respectively, and play an important role in
projective
differential geometry.  The vanishing of the Darboux cubic form
is characteristic for quadrics: indeed, in this case $\beta = \gamma =0$ so
that
asymptotic curves of both families are straight lines. The vanishing of the
projective metric (which is equivalent to either $\beta =0$ or $\gamma =0$)
characterises ruled surfaces. In what follows we exclude these two
degenerate
situations and require $\beta \ne 0$, $\gamma \ne 0$.

Using (\ref{newxy})-(\ref{new}), one can  verify that the four
points
\begin{equation}
\begin{array}{c}
{\bf r}, ~~~
{\bf r}_1={\bf r}_x-\frac{1}{2}\frac{\gamma_x}{\gamma}{\bf r}, ~~~
{\bf r}_2={\bf r}_y-\frac{1}{2}\frac{\beta_y}{\beta}{\bf r}, \\
\ \\
{\mbox{\boldmath $\eta$}}={\bf
r}_{xy}-\frac{1}{2}\frac{\gamma_x}{\gamma}{\bf
r}_y-
\frac{1}{2}\frac{\beta_y}{\beta}{\bf r}_x+
\left(\frac{1}{4}\frac{\beta_y\gamma_x}{\beta \gamma} -
\frac{1}{2}{\beta \gamma}\right){\bf r} \\
\end{array}
\label{frame}
\end{equation}
are defined in an invariant way, that is, under the transformation formulae
\mbox{(\ref{newxy})-(\ref{new})}
they acquire a nonzero multiple which does not
change them as points in projective space $P^3$. These points form the
vertices of the so-called Wilczynski moving tetrahedral \cite{Bol},
\cite{Finikov37}, \cite{Wilczynski}.
Since the lines passing through ${\bf r}, {\bf r}_1$ and ${\bf r}, {\bf
r}_2$
are
tangential to the $x$- and $y$-asymptotic curves,
respectively, the three points ${\bf r}, {\bf r}_1, {\bf r}_2$
span the tangent plane of the surface $M^2$.
The line through ${\bf r}_1, {\bf r}_2$ lying in the tangent
plane is
known as the directrix of Wilczynski of the second kind. The line through
${\bf r}, {\mbox{\boldmath $\eta$}}$ is transversal to $M^2$ and is known as
the directrix of
Wilczynski of the first kind. It plays the role of a projective `normal'.
 The Wilczynski tetrahedral proves to
be the most convenient tool in projective differential geometry.

Using (\ref{r}) and~(\ref{frame}),
we easily derive for ${\bf r}, {\bf r}_1, {\bf r}_2, {\mbox{\boldmath
$\eta$}}$
the linear equations \cite[p.\ 42]{Finikov37}
\begin{equation}
\begin{array}{c}
\left(\begin{array}{c}
{\bf r}\\
{\bf r}_1\\
{\bf r}_2\\
{\mbox{\boldmath $\eta$}}
\end{array}\right)_x=
\left(\begin{array}{cccc}
\frac{1}{2}\frac{\gamma_x}{\gamma}& 1&0&0\\
\frac{1}{2}b & -\frac{1}{2}\frac{\gamma_x}{\gamma}& \beta&0\\
\frac{1}{2}k&0&\frac{1}{2}\frac{\gamma_x}{\gamma}&1\\
\frac{1}{2}\beta
a&\frac{1}{2}k&\frac{1}{2}b&-\frac{1}{2}\frac{\gamma_x}{\gamma}
\end{array}\right)
\left(\begin{array}{c}
{\bf r}\\
{\bf r}_1\\
{\bf r}_2\\
{\mbox{\boldmath $\eta$}}
\end{array}\right)\\
\ \\
\left(\begin{array}{c}
{\bf r}\\
{\bf r}_1\\
{\bf r}_2\\
{\mbox{\boldmath $\eta$}}
\end{array}\right)_y=
\left(\begin{array}{cccc}
\frac{1}{2}\frac{\beta_y}{\beta}& 0&1&0\\
\frac{1}{2}l & \frac{1}{2}\frac{\beta_y}{\beta}& 0&1\\
\frac{1}{2}a&\gamma &-\frac{1}{2}\frac{\beta_y}{\beta}&0\\
\frac{1}{2}\gamma b&\frac{1}{2}a
&\frac{1}{2}l&-\frac{1}{2}\frac{\beta_y}{\beta}
\end{array}\right)\left(
\begin{array}{c}
{\bf r}\\
{\bf r}_1\\
{\bf r}_2\\
{\mbox{\boldmath $\eta$}}
\end{array}\right),
\end{array}
\label{Wilczynski}
\end{equation}
where we introduced the notation
\begin{equation}
\begin{array}{c}
k=\beta \gamma - (\ln \beta)_{xy}, ~~~~ l=\beta \gamma - (\ln
\gamma)_{xy},\\
\ \\
a=W-(\ln \beta)_{yy}-\frac{1}{2}(\ln \beta)_y^2, ~~~~
b=V-(\ln \gamma)_{xx}-\frac{1}{2}(\ln \gamma)_x^2.
\end{array}
\label{klab}
\end{equation}
The compatibility conditions of equations (\ref{Wilczynski}) imply
\begin{equation}
\begin{array}{c}
(\ln \beta)_{xy}=\beta \gamma -k, ~~~~ (\ln \gamma)_{xy}=\beta \gamma -l,\\
\ \\
a_x=k_y+\frac{\beta_y}{\beta}k, ~~~~ b_y=l_x+\frac{\gamma_x}{\gamma}l,\\
\ \\
\beta a_y+2a\beta_y=\gamma b_x+2b\gamma_x,
\end{array}
\label{GC2}
\end{equation}
which is just the equivalent form of the projective `Gauss-Codazzi'
equations
(\ref{GC1}).

Equations (\ref{Wilczynski}) can be rewritten in the Pl\"ucker coordinates.
For a convenience of the reader we briefly recall this construction.
Let us consider a line $l$ in $P^3$ passing through the points ${\bf a}$ and
${\bf b}$  with the
homogeneous coordinates ${\bf a}=(a^0:a^1:a^2:a^3)$ and ${\bf
b}=(b^0:b^1:b^2:b^3)$. With
the line $l$ we associate a point ${\bf a}\wedge {\bf b}$ in projective
space
$P^5$ with
the homogeneous coordinates
$$
{\bf a}\wedge {\bf b}=(p_{01}:p_{02}:p_{03}:p_{23}:p_{31}:p_{12}),
$$
where
$$
p_{ij}=\det
\left(\begin{array}{cc}
a^i & a^j \\
b^i & b^j
\end{array}\right).
$$
The coordinates $p_{ij}$ satisfy the well-known quadratic Pl\"ucker relation
\begin{equation}
p_{01}\, p_{23}+p_{02}\, p_{31}+p_{03}\, p_{12}=0.
\label{quadric}
\end{equation}
Instead of ${\bf a}$ and ${\bf b}$ we may consider an arbitrary linear
combinations thereof
without changing ${\bf a}\wedge{\bf b}$ as a point in $P^5$.
Hence, we arrive at the well-defined
Pl\"ucker coorrespondence $l({\bf a},{\bf b})\to {\bf a}\wedge {\bf b}$
between lines in
$P^3$ and points on the Pl\"ucker quadric in $P^5$.
Pl\"ucker correspondence plays an important role in the projective
differential geometry of surfaces and often sheds some new light on those
properties of surfaces which are not `visible' in $P^3$ but acquire a
precise
geometric meaning only in $P^5$. Thus, let us consider a surface
$M^2\in P^3$ with the Wilczynski tetrahedral
${\bf r}, {\bf r}_1, {\bf r}_2, {\mbox{\boldmath $\eta$}} $
satisfying equations (\ref{Wilczynski}).
Since the two pairs of points ${\bf r}, {\bf r}_1$ and ${\bf r}, {\bf r}_2$
generate two lines in $P^3$ which are tangential
to the $x$- and $y$-asymptotic curves, respectively, the formulae
$$
{\cal U}={\bf r} \wedge {\bf r}_1,\quad {\cal V}={\bf r} \wedge {\bf r}_2
$$
define the images of these lines under the Pl\"ucker embedding. Hence, with
any
surface $M^2\in P^3$ there are canonically associated two surfaces
${\cal U}(x, y)$ and ${\cal V}(x, y)$ in $P^5$ lying on the
Pl\"ucker quadric (\ref{quadric}).
In view of the formulae
$$
{\cal U}_x=\beta \, {\cal V},\quad {\cal V}_y=\gamma \, {\cal U},
$$
we conclude that the line in $P^5$ passing through a pair of points
$({\cal U}, {\cal V})$ can also be generated by the pair of points
$({\cal U}, {\cal U}_x)$ (and hence is
tangential to the $x$-coordinate line on the surface ${\cal U}$) or by a
pair
of
points $({\cal V}, {\cal V}_y)$ (and hence is tangential to the
$y$-coordinate
line on the
surface ${\cal V}$).
Consequently, the surfaces ${\cal U}$ and ${\cal V}$ are two focal surfaces
of
the congruence
of straight
lines $({\cal U}, {\cal V})$ or, equivalently,
${\cal V}$ is the Laplace transform of ${\cal U}$ with
respect to $x$ and ${\cal U}$ is the
Laplace transform of ${\cal V}$ with respect to $y$.
We emphasize that the $x$- and
$y$-coordinate lines on the surfaces ${\cal U}$ and ${\cal V}$
are not asymptotic but
conjugate. Continuation of the Laplace sequence in both directions, that is
taking the $x$-transform of ${\cal V}$, the $y$-transform of ${\cal U}$,
etc.,
leads, in the
generic case, to an infinite Laplace sequence in $P^5$ known as the Godeaux
sequence of a surface $M^2$ \cite[p.\ 344]{Bol}. The surfaces of the Godeaux
sequence carry important geometric information about the surface $M^2$
itself.

Introducing
$$
\begin{array}{c}
{\cal A} ={\bf r}_2\wedge {\bf r}_1+{\bf r}\wedge {\mbox{\boldmath $\eta$}},
~~~
{\cal B} ={\bf r}_1\wedge {\bf r}_2+{\bf r}\wedge {\mbox{\boldmath $\eta$}},
\\
\ \\
{\cal P} = 2\, {\bf r}_2\wedge {\mbox{\boldmath $\eta$}}, ~~~
{\cal Q} = 2\, {\bf r}_1\wedge {\mbox{\boldmath $\eta$}},
\end{array}
$$
we arrive at the following equations for the Pl\"ucker coordinates:

\begin{equation}
\begin{array}{c}
\left(\begin{array}{c}
{\cal U}\\
{\cal A}\\
{\cal P}\\
{\cal V}\\
{\cal B}\\
{\cal Q}
\end{array}\right)_x=
\left(\begin{array}{cccccc}
0 & 0 & 0 & \beta & 0 & 0\\
k & 0 & 0 & 0 & 0 & 0\\
0 & k & 0 & -\beta a & 0 & 0\\
0 & 0 & 0 & \frac{\gamma_x}{\gamma} & 1 & 0\\
0 & 0 & 0 & b & 0 & 1\\
-\beta a & 0 & \beta & 0 & b &-\frac{\gamma_x}{\gamma}
\end{array}\right)
\left(\begin{array}{c}
{\cal U}\\
{\cal A}\\
{\cal P}\\
{\cal V}\\
{\cal B}\\
{\cal Q}
\end{array}\right)\\
\ \\
\left(\begin{array}{c}
{\cal U}\\
{\cal A}\\
{\cal P}\\
{\cal V}\\
{\cal B}\\
{\cal Q}
\end{array}\right)_y=
\left(\begin{array}{cccccc}
\frac{\beta_y}{\beta} & 1 & 0 & 0 & 0 & 0\\
a & 0 & 1 & 0 & 0 & 0\\
0 & a & -\frac{\beta_y}{\beta} & -\gamma b & 0 & \gamma\\
\gamma & 0 & 0 & 0 & 0 & 0\\
0 & 0 & 0 & l & 0 & 0\\
-\gamma b & 0 & 0 & 0 & l & 0
\end{array}\right)
\left(\begin{array}{c}
{\cal U}\\
{\cal A}\\
{\cal P}\\
{\cal V}\\
{\cal B}\\
{\cal Q}
\end{array}\right).
\end{array}
\label{UAPVBQ}
\end{equation}
Equations (\ref{UAPVBQ}) are consistent with the following table
of scalar products:
\begin{equation}
({\cal U}, {\cal P})=-1, ~~~ ({\cal A}, {\cal A})=1, ~~~
({\cal V}, {\cal Q})=1, ~~~ ({\cal B}, {\cal B})=-1,
\label{table}
\end{equation}
all other scalar products being equal to zero. This defines a scalar
product of the signature (3, 3) which is the same as that of the quadratic
form
(\ref{quadric}).

Different types of surfaces can be defined by imposing additional
constraints
on $\beta$, $\gamma$, $V$, $W$ (respectively, $\beta, \gamma, k, l, a, b$),
so that, in a sense, projective differential geometry is the
theory of (integrable) reductions of the underdetermined system (\ref{GC1})
(respectively, (\ref{GC2})).
Although the three linear systems (\ref{r}),
(\ref{Wilczynski}) and (\ref{UAPVBQ}) are in fact equivalent, some of them
prove
to be more suitable for studying particular classes of projective surfaces
--- see \cite{Fer8}, \cite{Fer002} for the further discussion.

\section{Appendix B: the Lie sphere frame}

Here we describe the construction of the so-called Lie sphere frame
canonically associated with a surface in Lie sphere geometry (see
\cite{Fer002}). 
Although the construction follows essentially that of Blaschke
\cite{Blaschke}, our final formulae prove to be more suitable for
the purposes of this paper.
Let $M^2$ be a surface in $E^3$ parametrized by  coordinates $R^1, R^2$
of curvature lines, with the radius-vector ${\bf r}$ and the unit normal
${\bf n}$  satisfying the Weingarten equations (\ref{Weingarten}).
Introducing the 6-vectors
$$
{\bf U}=\left \{\frac{1+{\bf r}^2-2w^1({\bf r}, {\bf n})}{2}, ~~
\frac{1-{\bf r}^2+2w^1({\bf r}, {\bf n})}{2}, ~~ {\bf r}-w^1{\bf n}, ~~ w^1
\right \}
$$
and
$$
{\bf V}=\left \{\frac{1+{\bf r}^2-2w^2({\bf r}, {\bf n})}{2}, ~~
\frac{1-{\bf r}^2+2w^2({\bf r}, {\bf n})}{2}, ~~ {\bf r}-w^2{\bf n}, ~~ w^2
\right \},
$$
we readily  verify that
\begin{equation}
({\bf U},\ {\bf U})=({\bf U},\ {\bf V})=({\bf V},\  {\bf V})=0
\label{scalar}
\end{equation}
where the scalar product of 6-vectors is defined
by the indefinite quadratic form (\ref{Liequadric}). In what follows we use
the
same notation $(~,~)$ for both the scalar product defined by
(\ref{Liequadric})
as well as
for the standard Euclidean scalar product in $E^3$; however, the dimension
of  vectors will clearly indicate which one has to be choosen.

A direct computation gives
\begin{equation}
\begin{array}{c}
\partial_1{\bf U}=\partial_1w^1 \left \{ -({\bf r}, {\bf n}), ~
({\bf r}, {\bf n}), ~ -{\bf n}, ~ 1\right \}\\
\ \\
\partial_2{\bf U}=\partial_2w^1 \left \{ -({\bf r}, {\bf n}), ~
({\bf r}, {\bf n}), ~ -{\bf n}, ~ 1\right \}+
\frac{w^2-w^1}{w^2}\left \{ (\partial_2{\bf r}, {\bf r}), ~
-(\partial_2{\bf r}, {\bf r}), ~ \partial_2{\bf r}, ~ 0\right \} \\
\ \\
\partial_1{\bf V}=\partial_1w^2 \left \{ -({\bf r}, {\bf n}), ~
({\bf r}, {\bf n}), ~ -{\bf n}, ~ 1\right \}+
\frac{w^1-w^2}{w^1}\left \{ (\partial_1{\bf r}, {\bf r}), ~
-(\partial_1{\bf r}, {\bf r}), ~ \partial_1{\bf r}, ~ 0\right \} \\
\ \\
\partial_2{\bf V}=\partial_2w^2 \left \{ -({\bf r}, {\bf n}), ~
({\bf r}, {\bf n}), ~ -{\bf n}, ~ 1\right \}
\end{array}
\label{dUV}
\end{equation}
implying
\begin{equation}
\begin{array}{c}
\partial_1{\bf U}=\frac{\partial_1w^1}{w^1-w^2}({\bf U}-{\bf V}) \\
\ \\
\partial_2{\bf V}=\frac{\partial_2w^2}{w^2-w^1}({\bf V}-{\bf U}). \\
\end{array}
\label{12UV}
\end{equation}
Differentiating (\ref{scalar}) and taking into account
 (\ref{dUV}) and (\ref{12UV}), we conclude that the only nonzero
scalar products among the vectors ${\bf U}, ~ {\bf V}, ~ \partial_1{\bf U},
~ 
\partial_2{\bf U}, ~
\partial_1{\bf V}, ~ \partial_2{\bf V}$ are the following:
$$
\begin{array}{c}
(\partial_2{\bf U}, ~ \partial_2{\bf U})=(w^1-w^2)^2\, G_{22} \\
\ \\
(\partial_1{\bf V}, ~ \partial_1{\bf V})=(w^1-w^2)^2\, G_{11} .
\end{array}
$$
Here $G_{11}=(\partial_1{\bf n}, \partial_1{\bf n}), ~
G_{22}=(\partial_2{\bf n}, \partial_2{\bf n})$
are the components of the third fundamental form of the surface $M^2$.
Differentiating the zero scalar products among
${\bf U}$,  ${\bf V}$, $\partial_1{\bf U}$,  $\partial_2{\bf U}$,
$\partial_1{\bf V}$,  $\partial_2{\bf V}$ and keeping in mind  (\ref{12UV}),
one can show
that the triple ${\bf U}, \partial_2{\bf U}, \partial_2^2{\bf U}$
is orthogonal to the triple
${\bf V}, \partial_1{\bf V}, \partial_1^2{\bf V}$.
In order to complete the vectors ${\bf U}$ and ${\bf V}$ to a
frame in $P^5$ with the simplest possible table of scalar products, we will
choose
appropriate combinations among the triples ${\bf U}, \partial_2{\bf U},
\partial_2^2{\bf U}$
and ${\bf V}, \partial_1{\bf V}, \partial_1^2{\bf V}$, separately.
Up to a certain normalization, the choice
described below coincides with that from \cite{Blaschke}.

Let us introduce the normalized vectors
\begin{equation}
{\cal U}=\frac{{\bf U}}{\sqrt {G_{22}} (w^2-w^1)}, ~~~~
{\cal V}=\frac{{\bf V}}{\sqrt {G_{11}} (w^1-w^2)}.
\label{calUV}
\end{equation}
This normalization is convenient for several reasons: first of all,
equations
(\ref{12UV}) reduce to the Dirac equation
\begin{equation}
\begin{array}{c}
\partial_1{\cal U}=p\, {\cal V} \\
\ \\
\partial_2{\cal V}=q\, {\cal U} \\
\end{array}
\label{12calUV}
\end{equation}
with the coefficients $p$ and $q$ given by
$$
p=\frac{\partial_1w^1}{w^1-w^2}\frac{\sqrt {G_{11}}}{\sqrt {G_{22}}}, ~~~~
q=\frac{\partial_2w^2}{w^2-w^1}\frac{\sqrt {G_{22}}}{\sqrt {G_{11}}}.
$$
It is important that both $p$ are $q$ are
Lie-invariant  (we emphasize that coefficients in (\ref{12UV})
are not Lie-invariant). The reparametrization of coordinates
$$
(R^1)^* =f(R^1), ~~~ (R^2)^* = g(R^2)
$$
induces the transformation of $p$ and $q$ as follows:
\begin{equation}
p^{*}=p g'/(f')^2, ~~~~ q^{*}=q f'/(g')^2,
\label{new}
\end{equation}
so that we can introduce the Lie-invariant metric
$$
-pq\,  dR^1dR^2
$$
and the Lie-invariant cubic form
\begin{equation}
p\, (dR^1)^3-q\, (dR^2)^3
\label{cubicpq}
\end{equation}
(notice that only the conformal class of the cubic form does make an
invariant sense).

There exist one more important property of the normalized vector
${\cal U}$ (resp., ${\cal V}$). It turns out
that the action of the Lie sphere group in $E^3$ induces linear
transformations
of the coordinates of ${\cal U}$ (resp., ${\cal V}$).
Since this linear action should necessarily preserve the Lie quadric
(\ref{Liequadric}), we arrive at the well-known isomorphism of the Lie
sphere
group
and $SO(4, 2)$.  Thus, the normalization (\ref{calUV})
linearises the action of the Lie
sphere group (see \cite{Fer001} for the details).

The only nonzero scalar products among  normalized
vectors
${\cal U}, ~ {\cal V}, ~ \partial_1{\cal U}, ~ \partial_2{\cal U}, ~
\partial_1{\cal V}, ~ \partial_2{\cal V}$ are the following:
$$
(\partial_2{\cal U}, ~ \partial_2{\cal U})=
(\partial_1{\cal V}, ~ \partial_1{\cal V})=1.
$$
Obviously, the normalized triples
${\cal U}, \partial_2{\cal U}, \partial_2^2{\cal U}$ and
${\cal V}, \partial_1{\cal V}, \partial_1^2{\cal V}$
remain mutually orthogonal. Let us introduce the following vectors
${\cal A}, {\cal P}$ from the first triple:
$$
{\cal A}=\partial_2{\cal U}-\frac{\partial_2p}{p}\, {\cal U}, ~~~~~
{\cal P}=\partial_2{\cal A}-a\, {\cal U}
$$
which we require to have the following nonzero scalar products:
$$
({\cal A}, {\cal A})=1, ~~~~ ({\cal U}, {\cal P})=-1.
$$
This uniquely specifies
$$
a=-\frac{1}{2}\, (\partial_2{\cal A}, \partial_2{\cal A}).
$$
Similarly, we can choose the vectors
$$
{\cal B}=\partial_1{\cal V}-\frac{\partial_1q}{q}\, {\cal V}, ~~~~~
{\cal Q}=\partial_1{\cal B}-b\, {\cal V}
$$
with the nonzero scalar products
$$
({\cal B}, {\cal B})=1, ~~~~ ({\cal V}, {\cal Q})=-1,
$$
which fixes
$$
b=-\frac{1}{2}\, (\partial_1{\cal B}, \partial_1{\cal B}).
$$
Vectors ${\cal U}, ~ {\cal A}, ~ {\cal P}$ and
${\cal V}, ~ {\cal B}, ~ {\cal Q}$ constitute the  Lie sphere frame
with the following simple table of scalar products
\begin{equation}
({\cal A}, {\cal A})=1, ~~~ ({\cal U}, {\cal P})=-1, ~~~
({\cal B}, {\cal B})=1, ~~~ ({\cal V}, {\cal Q})=-1,
\label{scal}
\end{equation}
all other scalar products are zero, which is of the desired signature (4,
2).

Equations of motion of the Lie sphere frame can be conveniently represented
in
 matrix form (\ref{lie1}) (see \cite{Fer001})
$$
\begin{array}{c}
\partial_1\left(\begin{array}{c}
{\cal U}\\
{\cal A}\\
{\cal P}\\
{\cal V}\\
{\cal B}\\
{\cal Q}
\end{array}\right)=
\left(\begin{array}{cccccc}
0 & 0 & 0 & p & 0 & 0\\
k & 0 & 0 & 0 & 0 & 0\\
0 & k & 0 & -pa & 0 & 0\\
0 & 0 & 0 & \frac{\partial_1q}{q} & 1 & 0\\
0 & 0 & 0 & b & 0 & 1\\
pa & 0 & -p & 0 & b &-\frac{\partial_1q}{q}
\end{array}\right)
\left(\begin{array}{c}
{\cal U}\\
{\cal A}\\
{\cal P}\\
{\cal V}\\
{\cal B}\\
{\cal Q}
\end{array}\right)\\
\ \\
\partial_2\left(\begin{array}{c}
{\cal U}\\
{\cal A}\\
{\cal P}\\
{\cal V}\\
{\cal B}\\
{\cal Q}
\end{array}\right)=
\left(\begin{array}{cccccc}
\frac{\partial_2p}{p} & 1 & 0 & 0 & 0 & 0\\
a & 0 & 1 & 0 & 0 & 0\\
0 & a & -\frac{\partial_2p}{p} & qb & 0 & -q\\
q & 0 & 0 & 0 & 0 & 0\\
0 & 0 & 0 & l & 0 & 0\\
-qb & 0 & 0 & 0 & l & 0
\end{array}\right)
\left(\begin{array}{c}
{\cal U}\\
{\cal A}\\
{\cal P}\\
{\cal V}\\
{\cal B}\\
{\cal Q}
\end{array}\right).
\end{array}
$$
The compatibility conditions of  (\ref{lie1})
produce the equations (\ref{LieGC})
$$
\begin{array}{c}
\partial_1\partial_2 \ln p=pq -k, ~~~~
\partial_1\partial_2 \ln q=pq -l,\\
\ \\
\partial_1a=\partial_2k+\frac{\partial_2p}{p}\, k, ~~~~
\partial_2b=\partial_1l+\frac{\partial_1q}{q}\, l,\\
\ \\
p\,  \partial_2a+2\, a\, \partial_2p +
q\,  \partial_1b + 2\, b\, \partial_1q = 0,
\end{array}
$$
which can be viewed as the "Gauss-Codazzi" equations in Lie sphere geometry.
Another (equivalent) form of equations (\ref{LieGC}) can be obtained by
introducing $V$ and $W$
$$
\begin{array}{c}
V= b+\partial_1^2 \ln q+\frac{1}{2}\, (\partial_1 \ln q)^2, \\
\ \\
W=a+\partial_2^2 \ln p+\frac{1}{2}\, (\partial_2 \ln p)^2,
\end{array}
$$
which, upon the substitution into (\ref{LieGC}), implies (\ref{VW}):
$$
\begin{array}{c}
\partial_2^3p-2\, W\, \partial_2p-p\, \partial_2W+
\partial_1^3q-2\, V\, \partial_1q-q\, \partial_1V=0, \\
\ \\
\partial_1W=2\, q\, \partial_2p+p\, \partial_2q \\
\ \\
\partial_2V=2\, p\, \partial_1q+q\, \partial_1p.
\end{array}
$$

\section{Appendix C: the scalar products in $\Lambda^2({\bf C}^4)$}

Let us consider a space ${\bf C}^4$ equiped with the pseudo-Hermitian scalar
product of 
the signature $(2, 2)$
$$
({\bf a}, {\bf b})=-a^0\bar b^3+a^1\bar b^2+a^2\bar b^1-a^3\bar b^0,
$$
and define the wedge product ${\bf a}\wedge {\bf b}\in \Lambda^2({\bf C}^4)$
of the vectors
${\bf a}=(a^0, a^1, a^2, a^3)$ and ${\bf b}=(b^0, b^1, b^2, b^3)$ by the
formula
$$
{\bf a}\wedge {\bf b}=(y^0, y^1, y^2, y^3, y^4, y^5)
$$
where 
$$
\begin{array}{c}
y^0=\frac{1}{2}(p_{02}-p_{31}), ~~~ y^1=\frac{1}{2}(p_{02}+p_{31}), ~~~
y^2=\frac{1}{2}(p_{03}+p_{12}), \\
\ \\
 ~~~ y^3=\frac{1}{2i}(p_{03}-p_{12}), ~~~
y^4=\frac{1}{2i}(p_{01}-p_{23}), ~~~ y^5=\frac{1}{2i}(p_{01}+p_{23}).
\end{array}
$$
( $p_{ij}=a^ib^j-a^jb^i$). Let
$$
{\bf A}\wedge {\bf B}=(Y^0, Y^1, Y^2, Y^3, Y^4, Y^5)
$$
be the wedge product of any  other two vectors
${\bf A}=(A^0, A^1, A^2, A^3)$ and ${\bf B}=(B^0, B^1, B^2, B^3)$.
The pseudo-Hermitian scalar product $(\ , \ )$ in ${\bf C}^4$ induces
the pseudo-Hermitian scalar product $(\ , \ )$
in $\Lambda^2({\bf C}^4)$ as follows:
$$
({\bf a}\wedge {\bf b}, \ {\bf A}\wedge {\bf B})=
\frac{1}{2} 
 (({\bf a}, {\bf B})({\bf b}, {\bf A})-({\bf a}, {\bf A})({\bf b}, {\bf
B})).
$$
In terms of  $y^i$ and $  Y^i$ this pseudo-Hermitian scalar product takes
the form 
$$
({\bf a}\wedge {\bf b}, \ {\bf A}\wedge {\bf B})=
-y^0\bar Y^0+y^1\bar Y^1+y^2\bar Y^2+
y^3\bar Y^3+y^4\bar Y^4-y^5\bar Y^5,
$$
which is of the signature $(4, 2)$. Let us also define the complex scalar
product
$\{\ , \ \}$
in $\Lambda^2({\bf C}^4)$
$$
\{{\bf a}\wedge {\bf b}, \ {\bf A}\wedge {\bf B}\}=
\frac{1}{2} 
 ({\bf a}\wedge {\bf b} \wedge {\bf A}\wedge {\bf B})
$$
(here the right hand side is undestood as a half of the determinant of the
$4\times 4$ matrix with the rows  ${\bf a}, \
{\bf b}, \ {\bf A}, \ {\bf B}$).
In terms of  $y^i$ and $ Y^i$ this scalar product takes the form
$$
\{{\bf a}\wedge {\bf b}, \ {\bf A}\wedge {\bf B}\}=
-y^0 Y^0+y^1 Y^1+y^2 Y^2+
y^3 Y^3+y^4 Y^4-y^5 Y^5.
$$
Clearly,
$$
(\xi, \ \psi)=
\{\xi, \ \bar \psi\}
$$
for any $\xi, \psi$ in $\Lambda^2({\bf C}^4)$.

\section{Acknowledgements}

This research was supported by the EPSRC grant GR/N30941. I would like to
thank
F. Burstall for drawing my attention to the twistor theory and L. Mason for
usefull references. I would also like to thank A. Pushnitski and A.
Veselov for the discussions on Schr\"odinger operators
with magnetic fields. I am particularly grateful to K.R. Khusnutdinova for
a detailed investigation of the surfaces of revolution
appearing in section 4.

\end{document}